\documentclass[12pt]{amsart}
\usepackage{amsfonts,amssymb,amsmath,amscd,amsthm} %AMS package
\usepackage[dvipdfmx]{graphicx}
\usepackage{stmaryrd}
\usepackage[dvipdfmx,colorlinks,allcolors=blue]{hyperref}
\newtheorem{thm}{Theorem}
\newtheorem{lem}[thm]{Lemma}
\newtheorem{prop}[thm]{Proposition}

\newtheorem{cor}[thm]{Corollary}

\newtheorem{conj}[thm]{Conjecture}
\theoremstyle{definition}
\newtheorem{dfn}[thm]{Definition}
\newtheorem{ex}[thm]{Example}
\newtheorem{rmk}[thm]{Remark}

\numberwithin{thm}{section}
\numberwithin{equation}{section}

\newcommand{\Proof}{\noindent {\it Proof}.\ \ }
\newcommand{\Hom}{\operatorname{Hom}}
\newcommand{\Ext}{\operatorname{Ext}}
\newcommand{\Spec}{\operatorname{Spec}}

\newcommand{\ob}{\operatorname{ob}}
\newcommand{\Hilb}{\operatorname{Hilb}}
\newcommand{\HF}{\operatorname{HF}}
\newcommand{\red}{\operatorname{red}}
\newcommand{\im}{\operatorname{im}}
\renewcommand{\div}{\operatorname{div}}
\newcommand{\Bs}{\operatorname{Bs}}

\newcommand{\Pic}{\operatorname{Pic}}

\newcommand{\Aut}{\operatorname{Aut}}
\newcommand{\Fix}{\operatorname{Fix}}

\newcommand{\codim}[1]{\operatorname{codim}}

\renewcommand{\labelenumi}{{\rm (\arabic{enumi})}}

\newcommand{\mapright}[1]{%
  \smash{\mathop{%
      \hbox to 1cm{\rightarrowfill}}\limits^{#1}}}
\newcommand{\mapleft}[1]{%
  \smash{\mathop{%
      \hbox to 1cm{\leftarrowfill}}\limits^{#1}}}
\newcommand{\mapdown}[1]{\Big\downarrow
  \llap{$\vcenter{\hbox{$\scriptstyle#1\,$}}$ }}

\title[Obstructions to deforming space curves...]{
  Obstructions to deforming space curves
  lying on a smooth cubic surface}
\author{Hirokazu Nasu}
\subjclass[2010]{Primary 14C05; Secondary 14D15, 14H50}
\keywords{Hilbert scheme, space curve, obstruction, cubic surface}
\address{
  Department of Mathematical Sciences,
  Tokai University,
  4-1-1 Kitakaname, Hiratsuka, 
  Kanagawa 259-1292, JAPAN}
\email{nasu@tokai-u.jp}

\begin{document}

\begin{abstract}
In this paper, we study the deformations of 
curves in the projective $3$-space $\mathbb P^3$ (space curves),
one of the most classically studied objects in algebraic geometry.
We prove a conjecture due to J.~O.~Kleppe
(in fact, a version modified by Ph.~Ellia) 
concerning maximal families of space curves 
lying on a smooth cubic surface,
assuming the quadratic normality of its general members.
We also give a sufficient condition for curves lying on a cubic surface
to be obstructed in $\mathbb P^3$ in terms of lines on the surface.
For the proofs, we use the Hilbert-flag scheme of $\mathbb P^3$
as a main tool and apply a recent result on primary obstructions
to deforming curves on a $3$-fold developed 
by S.~Mukai and the author.
\end{abstract}

\maketitle

\section{Introduction}
\label{sec:intro}

Space curves, i.e.,~curves embedded into $\mathbb P^3$,
are one of the most classically studied objects in algebraic geometry
(cf.~\cite{Halphen,Noether}).
Among all space curves, curves lying on a smooth cubic surface
were intensively studied by virtue of a beautiful geometry endowed 
with the surface.
For example, Mumford~\cite{Mumford} found 
an example of a generically non-reduced component of the Hilbert scheme,
%in the celebrated example of 
whose general point corresponds to a space curve lying 
on a smooth cubic surface.
This example was beautifully generalized by Kleppe 
in his systematic study \cite{Kleppe87}
on $3$-maximal families of space curves.
(See also e.g.~\cite{Kleppe81,Gruson-Peskine82,Ellia87,Dolcetti-Pareschi88,Floystad93,MDP3,Nasu1,Vakil06,Mukai-Nasu,Nasu4,Kleppe-Ottem15,Dan17,Kleppe17,Nasu5,Nasu6,Nasu7} for further studies related to Mumford's example.)
Let $H(d,g)^{sc}$ denote the Hilbert scheme of 
smooth connected curves in $\mathbb P^3$ of degree $d$ and genus $g$.
Let $W$ be an irreducible closed subset of $H(d,g)^{sc}$.
Then the least degree $s(C)$ of surfaces containing a 
general member $C$ of $W$ is a basic invariant of $W$ and denoted by $s(W)$.
In this paper, $W$ is called a {\em $s$-maximal family (or subset)}
for $s\in \mathbb Z$, if $s(W)=s$ and if $W$ is maximal with respect to $s$,
i.e.,~$s(V)>s(W)$ for any irreducible closed subset $V$
containing $W$ properly.
Every irreducible component $V$ of $H(d,g)^{sc}$ is
a $s(V)$-maximal family, but the converse is not true.
Let $W\subset H(d,g)^{sc}$ be a $3$-maximal family
and suppose that its general member $C$ lies on a smooth cubic surface.
Kleppe \cite{Kleppe87} showed that if $d>9$ then $\dim W=d+g+18$,
and if moreover $g \ge 3d-18$ and $H^1(\mathbb P^3,\mathcal I_C(3))=0$, 
then $W$ is a generically smooth component of $H(d,g)^{sc}$.
Here and later, $\mathcal I_C$ denotes the sheaf of ideals defining $C$ 
in $\mathbb P^3$ and $\mathcal I_C(n):
=\mathcal I_C \otimes_{\mathbb P^3} \mathcal O_{\mathbb P^3}(n)$
for $n\in \mathbb Z$.
Moreover, he originated the following conjecture,
but here it is presented by modifications 
proposed by Ellia~\cite{Ellia87}.

\begin{conj}[Kleppe (a version modified by Ellia)]
  \label{conj:Kleppe}
  Suppose that $d>9$, $g \ge 3d-18$ and $C$ is linearly normal.
  If $H^1(\mathbb P^3,\mathcal I_C(3))\ne 0$,
  %i.e.,~$H^1(\mathbb P^3,\mathcal I_C(1))=0$.
  then 
  \begin{enumerate}
    \item \label{item:maximality}
    $W$ is an irreducible component of $(H(d,g)^{sc})_{\red}$ and
    \item \label{item:non-reduced}
    $H(d,g)^{sc}$ is generically non-reduced along $W$.
  \end{enumerate}
  %(with the reduced scheme structure)
  %of dimension $d+g+18$.
\end{conj}
Thus the $3$-maximal families $W$ in Conjecture~\ref{conj:Kleppe} 
are expected to give rise to generically non-reduced components 
of $H(d,g)^{sc}$.
The conclusion \eqref{item:maximality} of this conjecture
is equivalent to that
\begin{equation}
  \label{eqn:dimension}
  \dim_{[C]} H(d,g)^{sc}=d+g+18.
\end{equation}
The conclusion \eqref{item:non-reduced} follows from \eqref{item:maximality} 
because if $d >9$, then
\begin{equation}
  \label{eqn:tangential codimension}
  h^0(C,N_{C/\mathbb P^3})=\dim W+h^1(\mathbb P^3,\mathcal I_C(3)),
\end{equation}
where $h^0(C,N_{C/\mathbb P^3})$ represents the tangential dimension 
of $H(d,g)^{sc}$ at the point $[C]$ corresponding to $C$.
Ellia pointed out that \eqref{item:maximality}
is false if we drop the assumption of the linear normality of $C$
by counterexample
(see also Dolcetti-Pareschi~\cite{Dolcetti-Pareschi88}
  for more counterexamples).
The condition that $g\ge 3d-18$ is also necessary 
for \eqref{item:maximality} by dimension reason
%$\overline W$ to be a component in the Hilbert scheme
(cf.~Remark~\ref{rmk:dimension reason}).
Conjecture~\ref{conj:Kleppe} is related to a problem of classifying
all irreducible components $V$ of $H(d,g)^{sc}$
%($=\bigsqcup_{s\ge 1} \bigcup_{s(V)=s} V$) 
with $s(V)=s$,
so far this problem has been solved for $s\le 2$ 
(cf.~\cite{Tannenbaum80}, see also \cite[Prop.~4.11]{Nasu1})
and a very few (but partial) results are obtained 
for $s\ge 4$ (cf.~\cite{Kleppe-Ottem15,Nasu5,Kleppe17}).
%Including the original work of Kleppe himself, 
Several papers, e.g.~\cite{Kleppe87,Ellia87,Nasu1,Kleppe-Ottem15,Kleppe17} 
contributed to Conjecture~\ref{conj:Kleppe}.
It is known that
if $g$ is sufficiently large, %compared to $d$, 
then the conjecture holds to be true
(see Remark~\ref{rmk:known ranges}).
Mumford's example appears
in a region of $(d,g)$-plane for which the conjecture is known to be true,
and attains the minimal degree and the minimal genus 
in the region ($(d,g)=(14,24)$).
The main purpose of this paper is to settle down this conjecture
assuming further that $C$ is quadratically normal.

\begin{thm}
  \label{thm:main1}
  Conjecture~\ref{conj:Kleppe} 
  is true if $C$ is quadratically normal,
  i.e.,~$H^1(\mathbb P^3,\mathcal I_C(2))=0$.
\end{thm}

Note that if $d > 9$ then
the $3$-maximal families $W \subset H(d,g)^{sc}$
are in one-to-one correspondence
with the $7$-tuples $(a;b_1,\dots,b_6)$ of integers
satisfying certain numerical conditions 
(see \eqref{conds:standard-prescribed-nef-big} 
  in \S\ref{subsec:3-maximal families}).
Then for every $W$ (and $C$) in Conjecture~\ref{conj:Kleppe},
we have either $b_6=1$ or $b_6=2$
%and $C$ is quadratically normal if and only if $b_6=2$
(cf.~Lemmas~\ref{lem:abnormality} and \ref{lem:quadratically normal}).
Theorem~\ref{thm:main1} shows that
Conjecture~\ref{conj:Kleppe} is always true if $b_6=2$.

\medskip

Another purpose of this paper is to give a sufficient condition 
for curves $C$ lying on a smooth cubic surface
to be obstructed in $\mathbb P^3$.
Here we say $C$ is {\em (un)obstructed} in $\mathbb P^3$
if the Hilbert scheme of $\mathbb P^3$ is (non)singular at $[C]$.
Let $S$ be a smooth cubic surface in $\mathbb P^3$ and
$C$ a smooth connected curve on $S$.
Then since $-K_S$ is ample, we can easily see
that $H^1(C,N_{C/S})=0$
%=H^1(C,-K_S\big{\vert}_S+K_C)=0$
by adjunction.
Then it follows from the exact sequence 
$0 \rightarrow N_{C/S}\rightarrow N_{C/\mathbb P^3} \rightarrow N_{S/\mathbb P^3}\big{\vert}_C \rightarrow 0$
that $H^1(C,N_{C/\mathbb P^3})
\simeq H^1(C,N_{S/\mathbb P^3}\big{\vert}_C)$ and hence
every obstruction to deforming $C$ in $\mathbb P^3$
is contained in $H^1(C,N_{S/\mathbb P^3}\big{\vert}_C)$
(cf.~Remark~\ref{rmk:obstructions to deforming curves on a cubic}).
Let $L$ denote the class in $\Pic S$ of the invertible sheaf
$\mathcal O_S(C)\otimes_S N_{S/\mathbb P^3}^{-1}$ on $S$.
Then we have isomorphisms
$H^1(C,N_{S/\mathbb P^3}\big{\vert}_C)\simeq H^2(S,-L)$
and $H^1(\mathbb P^3,\mathcal I_C(3)) \simeq H^1(S,-L)$
(cf.~\eqref{isom:abnormality} and \eqref{isom:obstruction space}).
It follows from a general theory that
the Hilbert-flag scheme of $\mathbb P^3$ is nonsingular at $(C,S)$
(cf.~Lemma~\ref{lem:fano-fano})
%of expected dimesion $d+g+18$
and the first projection $pr_1$: $(C,S)\mapsto [C]$
from the scheme is smooth at $(C,S)$ if $H^1(S,-L)=0$.
This implies that $C$ is unobstructed in $\mathbb P^3$
if $H^i(S,-L)=0$ for either $i=1$ or $i=2$ (cf.~\cite{Kleppe87}).
Otherwise it follows from the Serre duality and a vanishing theorem
that $L+K_S$ is effective and $L$ is not nef
(cf.~Lemma~\ref{lem:not nef}).

\begin{thm}
  \label{thm:main2}
  Suppose that $L+K_S\ge 0$ %and $L^2>0$. %and $L$ is not nef.
  and there exists a (-1)-curve (i.e.~a line) $E$ on $S$
  such that $m:=-L.E>0$.
  Then $C$ is obstructed in $\mathbb P^3$ if either
  \begin{enumerate}
    \item $m=1$, or
    \item $2 \le m \le 3$ and
    the restriction map
    \begin{equation}
      \label{map:restriction}
      \varrho: H^0(S,\Delta) \rightarrow H^0(E,\Delta\big{\vert}_E)
    \end{equation}
    is surjective, where $\Delta:=L+K_S-2mE$ is a divisor on $S$.
  \end{enumerate}
\end{thm}
Some special cases of Theorem~\ref{thm:main2} were 
also proved in \cite{Dolcetti-Pareschi88} (for $m=3$) 
and ~\cite{Nasu1} (for $m=1$) (cf.~Remark~\ref{rmk:obstructed}).

In the proof of Theorems~\ref{thm:main1} and \ref{thm:main2},
we apply a recent result in \cite{Nasu5,Nasu8} 
(cf.~Theorem~\ref{thm:obstruction})
and prove that a part of the first order deformations $\tilde C$ of $C$
in $\mathbb P^3$ do not lift to any deformations
$\tilde {\tilde C}$ of $C$ over $k[t]/(t^3)$, 
where $k$ is the ground field. 
(Then $H(d,g)^{sc}$ is singular at $[C]$.)
In the case where $h^2(S,-L)=1$ 
(and hence $h^1(C,N_{C/\mathbb P^3})=1$),
we are even able to determine
the dimension of $H(d,g)^{sc}$ at $[C]$ 
(cf.~Proposition~\ref{prop:codimension}).
It is not easy to determine
the dimension of the Hilbert scheme at a given singular point.
Nevertheless, Theorems~\ref{thm:obstruction} and Lemma~\ref{lem:maximality}
make this determination possible
by a help of a geometry of lines on cubic surfaces.

The organization of this paper is as follows. 
In \S\ref{subsec:del Pezzo} we recall 
basic results on linear systems on del Pezzo surfaces, 
and prove a vanishing theorem
(cf.~Lemma~\ref{lem:crutial}), which is crucial to our proof of 
Theorem~\ref{thm:main1}.
In \S\ref{subsec:3-maximal families} we get more specialized into 
cubic surfaces and recall a well known correspondence between curves 
on a smooth cubic surface and $7$-tuples of integers.
In \S\ref{subsec:flag and primary} and \S\ref{subsec:obstructed}
we recall some results on Hilbert-flag schemes
and primary obstructions to deforming subschemes.
We prove Theorems~\ref{thm:main1} and \ref{thm:main2}
in \S\ref{sec:proof} and give some examples in \S\ref{sec:examples}.
Throughout the paper, 
we work over an algebraically closed field $k$ of characteristic $0$.

\vskip 3mm

\paragraph{\bf Acknowledgment}

I thank Prof.~Jan Oddvar Kleppe for 
many helpful discussions, his warm encouragement,
and showing me his unpublished preprint \cite{Kleppe83}.
Due to his comment, I obtained Proposition~\ref{prop:codimension}.
This paper was written during my stay at the University of Oslo (UiO)
in 2019 as a visiting researcher.
I thank Prof.~Kristian Ranestad and Prof.~John Christian Ottem
for inspiring discussions during the stay
and also UiO for providing the facilities.
I thank the referee for giving helpful comments improving 
the readability and quality of this paper.
An early proof of Theorem~\ref{thm:main1} was obtained during
my stay in UiO for the period from April to September in 2008.
Throughout the period, I was financially supported
by the Research Council of Norway project no.~188588.
This work was supported in part by 
JSPS KAKENHI Grant Numbers JP17K05210 and JP20K03541.

\section{Preliminaries}
\label{sec:preliminaries}

\subsection{Linear systems on del Pezzo surfaces}
\label{subsec:del Pezzo}

In this section, we collect some results on
linear systems on del Pezzo surfaces.
We refer to e.g.~\cite{Manin,Nasu4} for the proofs.

A {\em del Pezzo surface} is a smooth projective surface $S$ 
with ample anticanonical divisor $-K_S$.
Let $S$ be a del Pezzo surface over $k$.
Since $k$ is algebraically closed, 
$S$ is isomorphic to $\mathbb P^1\times \mathbb P^1$
or a blow-up of $\mathbb P^2$ at $r$ points ($r<9$) in general positions, 
i.e., no three are on a line, no six are on a conic 
and any cubic containing eight points is smooth at each of them
(cf.~\cite[\S24]{Manin}).
The self-intersection number $K_S^2$ is called the {\em degree} of $S$ 
and denoted by $\deg S$.
We have $\deg S=8$ if $S\simeq \mathbb P^1\times \mathbb P^1$ 
and $\deg S=9-r$ otherwise.
The anticanonical linear system $|-K_S|$ on $S$
is base point free if $\deg S\ge 2$,
and very ample if and only if $\deg S\ge 3$.
A curve $C$ on $S$ is called {\em a line} if $K_S.C=-1$ and $C^2=-1$
and {\em a conic} if $K_S.C=-2$ and $C^2=0$.
Thus there exist no lines on $S$ if 
$S\simeq \mathbb P^2$ or $\mathbb P^1\times \mathbb P^1$.
%To simplify our argument, in what follows, we assume 
%that $S\not\simeq \mathbb P^1\times \mathbb P^1$.

We recall some properties of divisors on a del Pezzo surface.
Let $D$ be a divisor on $S$.
We say that $D$ is nef if $D.C\ge 0$ for all curves $C$ on $S$.
Then we have the following:
\begin{enumerate}
  \item $D$ is nef if and only if $D \ge 0$
  and $D.\ell \ge 0$ for all lines $\ell$ on $S$\footnote{
    For $S$ isomorphic to neither
    $\mathbb P^2$ nor $\mathbb P^1\times \mathbb P^1$,
    $D$ is nef if and only if 
    $D.\ell \ge 0$ for all lines $\ell$ on $S$.
  }.
  \item If $D$ is nef, then $D^2 \ge 0$, where
  the equality holds if and only if there exists
  a conic $q$ on $S$ and an integer $m \ge 0$ such that $D \sim mq$.
  Then we say that $D$ is {\em composed with pencils}.
\end{enumerate}

Let $\chi(S,D)$ denote the Euler characteristic of
the invertible sheaf $\mathcal O_S(D)$ associated to $D$.
The following results are well-known.

\begin{lem}
  \label{lem:vanishing and base component}
  Let $D$ be a divisor on a del Pezzo surface $S$.
  Then
  \begin{enumerate}
    \item ({\em Zariski decomposition})
    Suppose that $D \ge 0$.
    Then for the complete linear system $|D|$ on $S$,
    %(spanned by $D$)
    there exists a unique decomposition
    $$
    |D| = |D'| + F,
    $$
    where $F$ is the fixed part of $|D|$
    (then $F$ is a $1$-dimensional subscheme of $S$).
    Here $D'$ is nef and $F$ is given by
    \begin{equation}
      \label{eqn:fixed part}
      F=-\sum_{D.\ell<0} (D.\ell) \ell,
    \end{equation}
    where the sum is taken over all lines $\ell$ on $S$ 
    such that $D.\ell<0$\footnote{
      The lines $\ell$ are mutually disjoint,
      i.e., if $(D.\ell)<0$ and $(D.\ell')<0$,
      then $\ell \cap \ell'=\emptyset$ .
      Thus the number of lines contained in the support of $F$ is at most $r$.
    }.
    In particular, $|D|$ is base point free if and only if $D$ is nef,
    except for the case where $\deg S=1$ and $D\sim -K_S$.
    \item Suppose that $D \ge 0$ and $D^2>0$.
    Then $H^1(S,-D)=0$ if and only if $D$ is nef.
    Otherwise, we have $h^1(S,-D)=h^0(F,\mathcal O_F)$,
    where $F=\Fix |D|$.
    \item If $D$ is nef and $\chi(S,-D)\ge 0$,
    then $H^1(S,-D)=0$.
  \end{enumerate}
\end{lem}
\Proof
(1) follows from \cite[Lemma~2.2]{Nasu4},
(2) from \cite[Lemma~2.4]{Nasu1}
and (3) from \cite[Lemma 2.1]{Nasu4}.
\qed

We apply the above lemma to the linear systems of
the canonical adjunctions and obtain the following corollary.

\begin{cor}[{cf.~\cite[Corollary~2.5]{Nasu1}}]
\label{cor:vanishing and base component}
Let $n$ be an integer and suppose that $D+nK_S\ge 0$.
Let $F$ be the fixed part of $|D+nK_S|$. Then
\begin{enumerate}
  \item $F= \sum_{D.\ell<n} (n-(D.\ell)) \ell$.
  \item Suppose that $(D+nK_S)^2>0$. 
  Then $H^1(S,-D-nK_S)=0$ if and only if $D+nK_S$ is nef.
  Otherwise, $h^1(S,-D-nK_S)=h^0(F,\mathcal O_F)$.
\end{enumerate}
\end{cor}
We will consider in \S\ref{sec:proof} the problem of determining 
the obstructedness of space curves for curves $C \subset \mathbb P^3$
lying on a smooth cubic surface $S \subset \mathbb P^3$.
Due to the following lemma,
we will be able to restrict ourselves
to the case where 
$\mathcal O_S(C)\otimes_S N_{S/\mathbb P^3}^{-1}$ is not nef 
(cf.~Theorem~\ref{thm:main2}).
\begin{lem}
  \label{lem:not nef}
  Let $D$ be a divisor on a del Pezzo surface $S$.
  If $H^i(S,-D)\ne 0$ for both $i=1$ and $i=2$,
  then $D+K_S$ is effective and $D$ is not nef.
\end{lem}
\Proof
By Serre duality, we have 
$H^0(S,D+K_S)\simeq H^2(S,-D)^{\vee}$,
which implies that $D+K_S\ge 0$.
Suppose that $D$ is nef for a contradiction.
Then $D^2\ge 0$. Since $H^1(S,-D)\ne 0$, this implies that
$D^2=0$ by Lemma~\ref{lem:vanishing and base component}.
Then $D$ is composed with pencils, i.e.,~there exists 
a conic $q$ on $S$ and an integer $m$
such that $D\sim mq$.
Since $q$ is nef, we see that $0\le q.(D+K_S)=mq^2+q.K_S=-2$,
thus a contradiction.
\qed

\medskip

We next setup some notations concerning the
coordinates of divisors on a del Pezzo surface.
Let $S$ be a blow-up of $\mathbb P^2$ at $r$ points in general position.
%which do not lie on a conic, and no three of which lie on a line.
Then the Picard group $\Pic S$ of $S$ has a $\mathbb Z$-free basis
$\mathbf l, \mathbf e_1, \dots, \mathbf e_r$
and we have $\Pic S\simeq \mathbb Z^{r+1}$.
Here and later, $\mathbf l$ and $\mathbf e_i$ $(1 \le i \le r)$ 
represent the class of the pullback of lines in $\mathbb P^2$
and $r$ exceptional curves on $S$, respectively.
Thus every divisor $D \sim a \mathbf l - \sum_{i=1}^r b_i \mathbf e_i$ on $S$
corresponds to a $(r+1)$-tuple $(a;b_1,\ldots,b_r)$ of integers by coordinates.
For examples, the anticanonical class 
$-K_S$ ($\simeq 3 \mathbf l - \sum_{i=1}^r \mathbf e_i$)
corresponds to $(3;1,\dots,1)$.
\begin{lem}[{cf.~\cite[V,Theorem~4.9]{Hartshorne}}] 
  \label{lem:lines and nefness}
  Suppose that $\deg S\ge 3$ and 
  $D \sim a \mathbf l - \sum_{i=1}^r b_i \mathbf e_i$ in $\Pic S$.
  \begin{enumerate}
    \item The class of lines on $S$ are represented by
    [i] $\mathbf e_i$ ($1 \le i \le r$),
    [ii] $\mathbf l -\mathbf e_i -\mathbf e_j$ for $1 \le i < j \le r$ and
    [iii] $2\mathbf l -\sum_{k=1}^5 \mathbf e_{i_k}$
    for $\left\{i_1,\dots,i_5\right\} \subset 
    \left\{1,\dots,r\right\}$.
    \item $D$ is nef if and only if
    \begin{enumerate}
      \item $b_i \ge 0$ for any integer $1 \le i \le r$,
      \item $a-b_i-b_j \ge 0$ for any integers $1 \le i < j \le r$ and
      \item $2a-\sum_{k=1}^5 b_{i_k}\ge 0$
      for any subset $\left\{i_1,\dots,i_5\right\} \subset 
      \left\{1,\dots,r\right\}$.
    \end{enumerate}
  \end{enumerate}
\end{lem}
We next take actions of Weyl groups into account.
For each $r\ge 2$, there exists a Weyl group $W_r \subset \Aut(\Pic S)$.
Here $W_r$ is generated by the permutations of $\mathbf e_i$ ($1 \le i \le r$)
and by the Cremona transformation $\sigma$ on $\mathbb P^2$
(only for $r\ge 3$),
where $\sigma$ is defined by
$\sigma(\mathbf l)=2\mathbf l-\mathbf e_1-\mathbf e_2-\mathbf e_3$
and 
$\sigma(\mathbf e_i)=\mathbf l-\sum_{1\le j\le 3, j\ne i}\mathbf e_j$
if $i \in \left\{1,2,3\right\}$ and
$\sigma(\mathbf e_i)=\mathbf e_i$ otherwise.
Then by virtue of this action of $W_r$ on $\Pic S$,
there exists a suitable blow-up $S \rightarrow \mathbb P^2$
such that
\begin{equation}\label{eqn:standard}
  b_1 \ge \dots \ge b_r \quad \mbox{and} \quad a \ge b_1 + b_2 + b_3
\end{equation}
(see e.g.~(cf.~\cite{Manin}) and \cite[\S5.3]{Nasu4}).
We say that the basis $\mathbf l,\mathbf e_1,\ldots,\mathbf e_r$
is {\em standard} for $D$ if \eqref{eqn:standard} is satisfied.
Under this basis, it is easily seen that 
$D$ is nef if and only if $b_r\ge 0$.
If $D$ is nef, then $|D|$ contains an irreducible curve as a member,
whose degree $d$ and (arithmetic) genus $g$ are respectively obtained 
by the formulas
\begin{equation}\label{eqn:degree and genus}
d = 3a- \sum_{i=1}^r b_i \quad \mbox{and} \quad 
g = \dfrac{(a-1)(a-2)}2 - \sum_{i=1}^r \dfrac{b_i(b_i-1)}2.
\end{equation}
Here the latter formula follows from the adjunction formula $2g-2=D.(K_S+D)$.
In particular, if $g>0$ then we have $a-b_1\ge 2$.
On the other hand, by the Hodge index theorem 
(cf.~\cite[Chap.~V,~Theorem~1.9]{Hartshorne}),
we have $(-K_S)^2D^2-(-K_S.D)^2=(9-r)(d+2g-2)-d^2\le 0$,
which implies 
\begin{equation}
  \label{ineq:Hodge index}
  0 \le g \le 1+(d-3)d/2
\end{equation}
if $r=6$ (i.e, $\deg S=3$). 

We need Lemma~\ref{lem:nef} below to prove Lemma~\ref{lem:crutial}.
For simplicity, we assume $\deg S\ge 3$ in these lemmas.
In what follows, the genus $g(D)$ of a divisor $D$ on $S$
(or an invertible sheaf $L$ on $S$) is defined by the adjunction formula,
which agrees with the genus of a curve $D$ if $D\ge 0$
(or $L\simeq \mathcal O_S(D)$).
\begin{lem}
  \label{lem:nef}
  Let $S$ be a del Pezzo surface,
  $F$ a divisor on $S$ that is a sum of mutually disjoint (single) $k$
  lines on $S$.
  Let $\varepsilon:S \rightarrow S'$ be the blow-down of $F$ from $S$
  and $D'$ a divisor on $S'$ of genus $g(D')\ge k$.
  If $D'$ is nef and $\deg S\ge 3$, then $\varepsilon^*D'-2F$ is nef.
  %  except for the following cases:
  %  \begin{enumerate}
    %    \item $|D|$ consists of a line on $S$ and $\deg S\le 2$, or
    %    \item $|D|=\emptyset$ and $\deg S=1$.
    %  \end{enumerate}
\end{lem}
\Proof
Let $E_k$ ($1 \le i \le k$) be mutually disjoint lines on $S$
and put $F:=\sum_{i=1}^k E_i$ and $D:=\varepsilon^*D'-2F$.
Then we note that $g(D)=g(D')-k$.
Therefore, it suffices to prove the lemma for $k=1$ by induction.
Suppose now that $F$ is a line on $S$ and $g(D')\ge 1$.
We put $r:=9-K_S^2$, i.e.,~the number of points of $\mathbb P^2$
blown-up to obtain $S$. Then since $\deg S\ge 3$, we have $1 \le r \le 6$.
By virtue of the action of the Weyl group $W_r$ on $\Pic S$,
we may assume that $\mathbf e_r$ is the class of $F$
and moreover, $D$ is linearly equivalent to
$a\mathbf l-\sum_{i=1}^{r-1} b_i\mathbf e_i -2\mathbf e_r$
with $a\ge b_1+b_2+b_3$ (only for $r\ge 4$)
and $b_1\ge \dots \ge b_{r-1}$.
Since $D'$ is nef, we have $b_{r-1} \ge 0$,
which implies $D.\mathbf e_i \ge 0$ for all $i$.
It follows from $g(D')\ge 1$ that 
$D.(\mathbf l-\mathbf e_i-\mathbf e_r)=a-b_i-2\ge a-b_1-2\ge 0$
for all $1 \le i\le r-1$.
Then we also see that 
$D.(2\mathbf l-\mathbf e_{i_1}-\dots-\mathbf e_{i_4}-\mathbf e_r)
=2a-b_{i_1}-\dots-b_{i_4}-2\ge2a-b_1-\dots-b_4-2\ge 0$
for all $1\le i_1 < \dots < i_4 \le r-1$.
Thus we have proved the lemma by Lemma~\ref{lem:lines and nefness}.
%Suppose that $K_S^2\le 2$.
%Then there exist also lines on $S$ corresponding to the classes
%$3\mathbf l-2\mathbf e_1-\sum_{i=2}^7 \mathbf e_i$ and
%$3\mathbf l-\sum_{i=1}^6 \mathbf e_i-2\mathbf e_r$.
%By considering the intersection number of $D$ with these lines,
%for proving that $D$ is nef, it is sufficient to show
%$$
%3a-2b_1-\dots-b_6-2\ge 0
%\qquad
%\mbox{and}
%\qquad
%3a-b_1-\dots-b_6-4\ge 0.
%$$
%The former inequality is alway satisfied, 
%while the latter is satisfied unless
%$a=3$ and $b_i=1$ for all $1 \le i \le 6$ (and with $b_7=0$ or $1$ for $r=8$).
%Then $D$ with such a tuple $(a;b_1,\dots,b_{r-1},2)$
%exactly corresponds to the exceptional cases mentioned in the lemma.
\qed

\begin{rmk}
  The conclusion of Lemma~\ref{lem:nef} is not true if $\deg S < 3$.
  In fact, suppose that $\deg S=2$ and $F$ is a line on $S$.
  Then $S'$ is a cubic surface.
  We consider the anticanonical class $D'=-K_{S'}$ on $S'$,
  whose genus is equal to one.
  Since $D=\varepsilon^* D'-2F$
  represents the class of a line on $S$, $D$ is clearly not nef.
\end{rmk}

The following lemma is a generalization of
Lemma~\ref{lem:vanishing and base component}
and plays an important role in our proof of Theorem~\ref{thm:main1}.

\begin{lem}
  \label{lem:crutial}
  Suppose that $\deg S\ge 3$ and $D$ is effective. If
  \begin{enumerate}
    \item $\chi(S,-D)\ge 0$ and
    %\item $D+K_S\ge 0$ (and hence $D \ge 0$) and
    \item $D.E \ge -1$ for any line $E$ on $S$,
  \end{enumerate}
  then we have $H^1(S,3F-D)=0$, where $F$ is the fixed part of $|D|$.
\end{lem}
\Proof If $D$ is nef, then we have $F=0$ and hence the lemma follows from
the first condition and Lemma~\ref{lem:vanishing and base component}.
Suppose now that $D$ is not nef.
Then by the same lemma, $F$ is non-empty.
It follows from the second condition that
$F$ is a sum $E_1+\dots+E_k$ of mutually disjoint lines $E_i$
($1 \le i \le k$) on $S$.
Let $\varepsilon:S\rightarrow S'$ be the blow-down of $F$ from $S$.
Then $\Delta:=D-F$ is the pull-back of a nef divisor $\Delta'$ on $S'$, 
i.e.~$\varepsilon^*\Delta'\sim \Delta$.
Since $D\sim \varepsilon^*\Delta'+F$, $K_S=\varepsilon^*K_{S'}+F$ and $F^2=-k$,
%we have $D^2=\Delta'^2-k$ and $D.K_S=\Delta'K_{S'}-k$. Then 
we compute that
$$
2\chi(S,-D)-2=D.(D+K_S)=\Delta'.(\Delta'+K_{S'})-2k=2\chi(S',-\Delta')-2-2k.
$$
Thus $g(\Delta')=\chi(S',-\Delta')=\chi(S,-D)+k\ge k$.
Moreover, since $\Delta'$ is nef,
so is $D-3F\sim \varepsilon^*\Delta'-2F$ by Lemma~\ref{lem:nef}.
Then we note that $D.F=F^2=-k$
and hence we see that $\chi(S,3F-D)\ge 0$ by
$$
\chi(S,3F-D)-\chi(S,-D)=\chi(S,3F)-\chi(S,\mathcal O_S)-3D.F=-3k+3k=0.
$$
Thus we have completed the proof
by Lemma~\ref{lem:vanishing and base component}.
\qed

\medskip

Finally, we prepare a lemma concerning the bigness of divisors.

\begin{lem}
  \label{lem:big}
  If $\deg S\ge 3$, $D+K_S\ge 0$ and $\chi(S,-D)\ge 0$
  %and $D$ is not nef 
  then $D$ is big, i.e.,~$D^2>0$.
\end{lem}
\Proof
It follows from the Riemann-Roch theorem that
$$
D^2=2\chi(S,-D)-2-D.K_S\ge -2-K_S.D.
$$
Then since $-K_S$ is ample, we have
$$
-2-K_S.D=-2-K_S.(D+K_S)+K_S^2\ge -2+K_S^2>0. \qed
$$
%Then since $D$ is not nef,
%there exists a line $E$ on $S$ such that $D.E<0$.
%Since $(D+K_S.E)=D.E-1<-1$, we see that
%$D+K_S-2E$ is effective (cf.~Corollary~\ref{cor:vanishing and base component}).
%Then since $-K_S$ is ample, we have
%$$
%-2-K_S.D=-2-K_S(D+K_S-2E)+K_S^2+2\ge K_S^2>0.
%$$

\subsection{3-maximal families}
\label{subsec:3-maximal families}

In this section, we consider cubic del Pezzo surfaces.
First we recall some properties of curves on the surfaces.
Let $S$ be a smooth cubic surface.
Then $S$ is a blow-up of $\mathbb P^2$ at six points in general position.
As we have seen in the previous section,
for every divisor $D$ on $S$,
$\Pic S$ has a standard basis $\mathbf l,\mathbf e_1,\ldots,\mathbf e_6$
and $D$ corresponds to
a $7$-tuple $(a;b_1,\ldots,b_6)$ of integers
satisfying
\begin{equation}\label{eqn:standard2}
  b_1 \ge \dots \ge b_6 \quad \mbox{and} \quad a \ge b_1 + b_2 + b_3.
\end{equation}
Then $D$ is nef if and only if $b_6\ge 0$.
Moreover, $|D|$ contains a smooth connected curve not a line nor a conic
if and only if $a > b_1$ and $b_6 \ge 0$.
Its degree $d$ and genus $g$ 
are respectively given by \eqref{eqn:degree and genus}
and they satisfy the inequality \eqref{ineq:Hodge index}.

We next consider the projective (ab)normality of curves 
on a smooth cubic surface.
Given an integer $n$, 
a projective variety $V \subset \mathbb P^d$
is said to be {\em $n$-normal} if $H^1(\mathbb P^d,\mathcal I_V(n))=0$.
$V$ is called {\em projectively normal}
if $V$ is $n$-normal for all $n \in \mathbb Z$.

\begin{lem}
  \label{lem:abnormality}
  Let $C$ be a curve on a smooth cubic surface $S$.
  We assume that $C+nK_S\ge 0$ and $(C+nK_S)^2>0$.
  Then
  $$
  h^1(\mathbb P^3,\mathcal I_C(n))=h^0(F,\mathcal O_F),
  $$
  where $F$ is the fixed part of $|C+nK_S|$.
  In particular, $C$ is $n$-normal if and only if $C+nK_S$ is nef.
\end{lem}

\medskip

\Proof 
We note that $\mathcal I_S(n)\simeq \mathcal O_{\mathbb P^3}(n-3)$,
whose $i$-th cohomology groups vanish for $i=1,2$.
It follows from an exact sequence 
$0 \rightarrow \mathcal I_S(n)
  \rightarrow \mathcal I_C(n)
  \rightarrow \mathcal O_S(-nK_S-C)
  \rightarrow 0
  $
  on $\mathbb P^3$ that
  \begin{equation}
    \label{isom:abnormality}
    H^1(\mathbb P^3,\mathcal I_C(n))
    \simeq H^1(S,-nK_S-C).
  \end{equation}
  Therefore, the lemma follows 
  from Corollary~\ref{cor:vanishing and base component}.
\qed

\begin{rmk}
  \label{rmk:normality}
  Since $-K_S$ is effective and ample,
  if $D$ is nef and big, then so is $D-K_S$.
  This implies that in the setting of Lemma~\ref{lem:abnormality},
  if $C$ is $n$-normal then $C$ is $m$-normal for all $m < n$.
\end{rmk}

Let $C$ be a curve on $S$ with coordinate $(a;b_1,\dots,b_6)$ 
under a standard basis.
We note that every element of the Weyl group preserves the class $K_S$.
Therefore, if $\mathbf l,\mathbf e_1,\ldots,\mathbf e_6$ 
is a standard basis for $C$, then so is for $D:=C+nK_S$ with any integer $n$.
Therefore, provided that $D \ge 0$ and $D^2>0$,
$C$ is $n$-normal if and only if $b_6\ge n$ by Lemma~\ref{lem:abnormality}.
One should be careful in applying Lemma~\ref{lem:abnormality}
to a computation of the $n$-abnormality
$h^1(\mathbb P^3,\mathcal I_C(n))$ of $C$.
The support of the fixed part $F$ of $|C+nK_S|$
consists of any set of mutually disjoint lines,
whose number is at most $6$.
The next example shows that even under the standard basis,
the support of $F$ does not necessarily consist of
$\mathbf e_1,\dots,\mathbf e_6$.
This fact corresponds to the fact that 
every blow-down of a cubic surface along $5$ lines is isomorphic to
either $\mathbb P^1 \times \mathbb P^1$, or $\mathbb P^2$ blown up at a point.
%($\simeq \mathbb P(\mathcal O_{\mathbb P^1} \oplus \mathcal O_{\mathbb P^1}(-1))$

\begin{ex}
  \label{ex:fixed part}
  Let $C$ be a curve (of $d=18$ and $g=31$)
  on a smooth cubic surface $S$
  corresponding to the $7$-tuple $(12;5,5,2,2,2,2)$.
  Here and later, we abuse notations
  and identify divisor classes on $S$
  with $7$-tuples of integers corresponding to them.
  We put $D:=C+3K_S$ in $\Pic S$. Then we see that
  $$
  D=(2;1,1,0,0,0,0)+(1;1,1,-1,-1,-1,-1)=D'+F,
  $$
  where $D'=2\mathbf l-\mathbf e_1-\mathbf e_2$ is nef
  and $F$ is the fixed part of $|C+3K_S|$.
  We note that $F$ consists of $5$ disjoint lines
  $\mathbf l-\mathbf e_1-\mathbf e_2$ and $\mathbf e_i$ ($3\le i\le 6$).
  Since $D'$ is also big by $D'^2=2$, we have $H^1(S,-D')=0$.
  Then it follows from \eqref{isom:abnormality} and the exact sequence
  $
  0 \rightarrow \mathcal O_S(-D) \rightarrow \mathcal O_S(-D')
  \rightarrow \mathcal O_F \rightarrow 0
  $
  that
  $h^1(\mathbb P^3,\mathcal I_C(3))=h^1(S,-D)=h^0(F,\mathcal O_F)=5$.
\end{ex}

Let $d>0$ and $g$ be two integers satisfying \eqref{ineq:Hodge index}
and $(a;b_1,\ldots,b_6)$ a $7$-tuple of integers 
satisfying a set of conditions
\begin{equation}
  \label{conds:standard-prescribed-nef-big}
  \mbox{\eqref{eqn:standard2}, \eqref{eqn:degree and genus} (\mbox{with $r=6$}),
    $a > b_1$ and $b_6 \ge 0$}.
  %\tag{*}
\end{equation}
Then according to \cite{Kleppe87}, we can associate to $(a;b_1,\ldots,b_6)$
a closed subset of the Hilbert scheme.
Let $H(d,g)^{sc}$ be the Hilbert scheme of smooth connected curves
of degree $d$ and genus $g$ in $\mathbb P^3$.

\begin{dfn}\label{dfn:3-maximal}
  We define a closed subset $W(a;b_1,\dots,b_6) \subset H(d,g)^{sc}$
  by taking the closure in $H(d,g)^{sc}$
  of the family of curves $C \subset \mathbb P^3$ lying on a
  smooth cubic surface $S \subset \mathbb P^3$
  and such that 
  $$
  C \sim a \mathbf l - \sum_{i=1}^6 b_i \mathbf e_i
  $$
  on $S$ for some (standard) basis $\mathbf l, \mathbf e_1,\dots, \mathbf e_6$ of $\Pic S$.
\end{dfn}

Let $W=W(a,b_1,\dots,b_6)$.
If $d>9$ then every general member $C$ of $W$
is contained in a unique cubic surface $S$, and hence
$W$ is birationally equivalent to 
$\mathbb P^{d+g-1}$-bundle over 
$|\mathcal O_{\mathbb P^3}(3)|\simeq \mathbb P^{19}$, 
where the numbers $d+g-1$ and $19$ are equal to
the dimensions of the linear systems
$|\mathcal O_S(C)|$ on $S$ and 
$|\mathcal O_{\mathbb P^3}(3)|$ on $\mathbb P^3$, respectively.
In particular, $W$ is irreducible and of dimension $d+g+18$.
It is known that if $d>2$ then
every $3$-maximal family in $H(d,g)^{sc}$
(see ~\S\ref{sec:intro} for its definition) 
can be obtained as $W(a;b_1,\dots,b_6)$ for some $(a;b_1,\dots,b_6)$
satisfying \eqref{conds:standard-prescribed-nef-big},
provided that its general member is contained in a smooth cubic surface
(cf.~\cite{Kleppe87,Kleppe81}).
%because it is obtained as the projection 
%of some irreducible component of the Hilbert-flag scheme.
Conversely, if $d>9$ then $W(a;b_1,\ldots,b_6)$ becomes a $3$-maximal family.

\begin{rmk}
  \label{rmk:dimension reason}
  By deformation theory, every irreducible component of 
  $H(d,g)^{sc}$ is of dimension at least $4d$ 
  ($=\chi(C,N_{C/\mathbb P^3})$).
  Therefore, if $d >9$ and $W=W(a;b_1,\ldots,b_6) \subset H(d,g)^{sc}$
  is an irreducible component of $(H(d,g)^{sc})_{\red}$,
  then we have $g\ge 3d-18$ by dimension.
\end{rmk}

\subsection{Hilbert-flag schemes and Primary obstructions}
\label{subsec:flag and primary}

In this section, we briefly recall the definition of Hilbert-flag schemes
and their infinitesimal properties (cf.~\cite{Kleppe87,Sernesi}).
Given a projective scheme $X$ and a pair of Hilbert polynomials $P$ and $Q$,
there exists a contravariant functor
$HF_{P,Q}: (\mbox{schemes}) \rightarrow (\mbox{sets})$
that to each base scheme $B$ assigns a
pair of closed subschemes $\mathcal C \subset \mathcal S \subset X\times_k B$,
both flat over $B$, and where the fibers of $\mathcal C$ (resp., $\mathcal S$)
have the Hilbert polynomial $P$ (resp. $Q$).
This functor is represented by a projective scheme $\HF_{P,Q} X$,
so called the {\em Hilbert-flag scheme} of $X$.
Let $(C,S)$ be a pair of closed subschemes of $X$
with Hilbert polynomials $(P,Q)$, respectively and such that $C \subset S \subset X$.
Then the normal sheaf $N_{(C,S)/X}$ of $(C,S)$
is a sheaf of $\mathcal O_X$-module and defined by the fiber product
$$
N_{(C,S)/X} :=N_{C/X} \times_{N_{S/X}\big{\vert}_C} N_{S/X}
$$
of the projection $N_{C/X} \rightarrow N_{S/X}\big{\vert}_C$
and the restriction $N_{S/X} \rightarrow N_{S/X}\big{\vert}_C$
of normal sheaves of $C$ and $S$ in $X$, respectively 
(cf.~\cite[\S2.2]{Nasu6}).

In what follows,
we assume that the two embeddings $C \hookrightarrow S$ and 
$S \hookrightarrow X$ are both regular
(then so is $C \hookrightarrow X$).
Then it follows from a general theory 
(cf.~\cite[Proposition~4.5.3]{Sernesi})
that $H^0(X,N_{(C,S)/X})$ and $H^1(X,N_{(C,S)/X})$ respectively
represent the tangent space and the obstruction space
of $\HF_{P,Q} X$ at $(C,S)$, and we have
\begin{equation}
  \label{ineq:dimension of flag}
  h^0(X,N_{(C,S)/X})-h^1(X,N_{(C,S)/X}) \le \dim_{(C,S)} \HF_{P,Q} X
  \le h^0(X,N_{(C,S)/X}).
\end{equation}
Thus if $H^1(X,N_{(C,S)/X})=0$ then 
$\HF_{P,Q} X$ is nonsingular at $(C,S)$ of expected dimension,
that is, the number in the left hand side of \eqref{ineq:dimension of flag}.
The expected dimension coincides with $\chi(X,N_{(C,S)/X})$,
provided that $H^i(X,N_{(C,S)/X})=0$ for all $i \ge 2$.
Let $\Hilb_P X$ denote the Hilbert scheme of $X$ with Hilbert polynomial $P$.
Then there exist two natural projections
$pr_1: \HF_{P,Q} X \rightarrow \Hilb_P X$ and
$pr_2: \HF_{P,Q} X \rightarrow \Hilb_Q X$,
i.e.,~the first and the second projections.
Correspondingly, there exist two natural exact sequences
\begin{equation}
  \label{ses:normal1}
  \begin{CD}
      0 @>>> \mathcal I_{C/S}\otimes_S N_{S/X} 
      @>>> N_{(C,S)/X} @>{\pi_1}>> N_{C/X} @>>> 0
  \end{CD}
\end{equation}
and
\begin{equation}
  \label{ses:normal2}
  \begin{CD}
  0 @>>> N_{C/S} @>>> N_{(C,S)/X} @>{\pi_2}>> N_{S/X} @>>> 0
  \end{CD}
\end{equation}
of sheaves on $X$, where $\pi_1$ and $\pi_2$
induce the tangent map and the map on obstruction spaces of 
$pr_1$ and $pr_2$, respectively (cf.~\cite[\S2.2]{Nasu6}).

We recall that a normal projective variety $Z$ is called 
{\em Fano} if $-K_Z$ is ample. The following lemma shows that
if $S$ and $X$ are both Fano and if the two embeddings
$C\hookrightarrow S$ and $S\hookrightarrow X$ are both of codimension one,
then all the higher cohomology groups of $N_{(C,S)/X}$ vanish,
and we benefit a nice property from the Hilbert-flag scheme of $X$.

\begin{lem}
  \begin{enumerate}
    \item If $S$ and $X$ are both Fano,
    and both $C \subset S$ and $S \subset X$ 
    are effective Cartier divisors,
    then we have $H^i(X,N_{(C,S)/X})=0$ for all $i >0$.
    \item If $X$ is a Fano $3$-fold, $S$ is a del Pezzo surface and 
    $C$ is a curve of degree $d=(-K_S.C)$ and genus $g$,
    then $\HF X$ is nonsingular at $(C,S)$ of expected dimension
    $$
    \chi(X,N_{(C,S)/X})=\frac{(-K_X.S^2)_X}2+d+g,
    $$
    where $(D_1.D_2.D_3)_X$ denotes the intersection number of
    divisors $D_1,D_2,D_3$ on $X$.
  \end{enumerate}
  \label{lem:fano-fano}
\end{lem}
\Proof 
We note by adjunction that $N_{C/S}\simeq -K_S\big{\vert}_C+K_C$
and $N_{S/X}\simeq -K_X\big{\vert}_S+K_S$.
Therefore, the higher cohomology groups
$H^i(C,N_{C/S})$ and $H^i(S,N_{S/X})$ vanish for all $i>0$ 
by the ampleness of $-K_S$ and $-K_X$, respectively.
Thus (1) follows from the exact sequence \eqref{ses:normal2}.
By Riemann-Roch formulas on curves and surfaces, we see that
$\chi(C,N_{C/S})=d+g-1$
and $\chi(S,N_{S/X})=(-K_X.S^2)_X/2+1$.
Hence we obtain (2) by additivity on Euler characteristics.
\qed
%$\chi(C,N_{C/S})=\deg(-K_S\big{\vert}_C+K_C)=d+g-1$.
%$\chi(S,N_{S/X})=(O_S(S).O_S(S)-K_S)/2+1=(O_S(S).-K_X\big{\vert}_S)/2+1=(-K_X.S^2)_X/2+1

\begin{lem}[{cf.~\cite{Kleppe87,Nasu6}}]
  \label{lem:expected dimension}
  If $S \subset \mathbb P^3$ is a smooth cubic surface,
  and $C$ is a smooth curve on $S$ of degree $d$ and genus $g$, then
  \begin{enumerate}
    \item $H^1(C,N_{C/S})=H^1(S,N_{S/\mathbb P^3})=0$.
    \item $\HF \mathbb P^3$ is nonsingular at $(C,S)$ 
    of expected dimension
    $\chi(\mathbb P^3,N_{(C,S)/\mathbb P^3})=d+g+18$,
    which coincides with the dimension of
    $3$-maximal families in $H(d,g)^{sc}$ containing $C$
    if $d>9$ (cf.~\S\ref{subsec:3-maximal families}).
  \end{enumerate}
\end{lem}
\Proof
Since $S$ and $\mathbb P^3$ are both Fano, we obtain (1).
Since $-K_{\mathbb P^3}\sim 4H$ and $S \sim 3H$, where 
$\Pic \mathbb P^3 \simeq \mathbb Z [H]$, 
we see that $(-K_{\mathbb P^3}.S^2)_{\mathbb P^3}=36$ and 
thus (2) follows from Lemma~\ref{lem:fano-fano}. \qed

\medskip

We next recall the definition of primary obstructions
to deforming subschemes.
Let $\tilde C$ be a {\em first order deformation} of $C$ in $X$,
that is, an infinitesimal deformation $\tilde C$ of $C$ in $X$
over the ring $k[t]/(t^2)$ of dual numbers.
Then $\tilde C$ naturally corresponds to a global section 
$\alpha$ of $N_{C/X}$.
Since the embedding $C \hookrightarrow X$ is regular,
every obstruction to deforming $C$ in $X$ is contained in $H^1(C,N_{C/X})$
(cf.~\cite[Theorem~4.3.5]{Sernesi}).
Every $\alpha$ in $H^0(C,N_{C/X})$
defines an element $\ob(\alpha)$ of $H^1(C,N_{C/X})$
such that $\ob(\alpha)$ is zero if and only if
$\tilde C$ extends to a deformation 
${\tilde {\tilde C}}$ of $C$ over $k[t]/(t^3)$.
Here $\ob(\alpha)$ is called the {\em primary obstruction} of $\alpha$
(or $\tilde C$).
It is known that $\ob(\alpha)$ is expressed as 
a cup product of
$\alpha \in \Hom_X (\mathcal I_C,\mathcal O_C)\simeq H^0(C,N_{C/X})$
and the extension class 
$\mathbf e =[0 \rightarrow \mathcal I_C \rightarrow \mathcal O_X
\rightarrow \mathcal O_C \rightarrow 0]
\in \Ext^1(\mathcal O_C,\mathcal I_C)$
and we have $\ob(\alpha)=\alpha \cup \mathbf e \cup \alpha$
(cf.~\cite[Theorem~2.1]{Nasu5}).
If $\ob(\alpha)\ne 0$
then $\tilde C$ does not lift to a global deformation of $C$ in $X$
and $\Hilb_P X$ is singular at $[C]$.

\begin{rmk}
  \label{rmk:obstructions to deforming curves on a cubic}
  Here we give a remark on 
  obstructions to deforming space curves lying on a smooth cubic surface.
  Let $S \subset \mathbb P^3$ be a smooth cubic surface, 
  $C$ a smooth curve on $S$.
  Then by Lemma~\ref{lem:expected dimension}, 
  we always have $H^1(C,N_{C/S})=H^1(S,N_{S/\mathbb P^3})=0$.
  Thus both of the deformations of $C$ in $S$ and 
  the deformations of $S$ in $\mathbb P^3$
  behave well. However, those of $C$ in $\mathbb P^3$ can behave badly in general,
  For example, for curves $C$ considered in Mumford's example~\cite{Mumford},
  i.e., $C \sim -4K_S+2E$ and $E$ is a line on $S$, 
  we see that $H^1(C,N_{C/\mathbb P^3})\ne 0$.
  It follows from the exact sequence 
  $0 \rightarrow N_{C/S}\rightarrow N_{C/\mathbb P^3} \rightarrow 
  N_{S/\mathbb P^3}\big{\vert}_C \rightarrow 0$
  that $H^1(C,N_{C/\mathbb P^3})
  \simeq H^1(C,N_{S/\mathbb P^3}\big{\vert}_C)$ and hence
  every obstruction to deforming $C$ in $\mathbb P^3$
  is contained in $H^1(C,N_{S/\mathbb P^3}\big{\vert}_C)$,
  and this cohomology group does not vanish in the case of Mumford's example.
  In fact, it was proved by Curtin~\cite{Curtin} that
  there exists a first order deformation $\tilde C$ of $C$ in $\mathbb P^3$
  whose primary obstruction is nonzero in $H^1(C,N_{C/\mathbb P^3})$.
  Our method of computing primary obstructions,
  which will be explained in \S~\ref{subsec:obstructed}, is based on a technique
  used in \cite{Curtin} and also its generalization in \cite{Nasu1,Mukai-Nasu}.
\end{rmk}

Let $\mathcal W_{C,S}$ be an irreducible component of $\HF_{P,Q} X$
passing through $(C,S)$,
and let $pr_1': \mathcal W_{C,S} \rightarrow \Hilb_P X$ be 
the restriction of $pr_1$ to $\mathcal W_{C,S}$.
Then it follows from a general deformation theory 
(cf.~\cite[Lemma~A10]{Kleppe87}, see also \cite[Theorem~1.3.4]{Kleppe81})
that if $H^1(S,\mathcal I_{C/S}\otimes_S N_{S/X})=0$,
then $pr_1$ is smooth at $(C,S)$. Then so is $pr_1'$.
If moreover $H^1(X,N_{(C,S)/X})=0$
then $pr_1'$ is dominant in a neighborhood of $[C]$
(cf.~\cite[Theorem~2.4]{Nasu6}).
If $H^1(X,N_{(C,S)/X})=0$ and 
$H^1(S,\mathcal I_{C/S}\otimes_S N_{S/X})\ne 0$,
then there exists an exact sequence
\begin{equation}
  \label{ses:first tangent map}
  \begin{CD}
    H^0(X,N_{(C,S)/X})
    @>{p_1}>> H^0(C,N_{C/X})
    @>{\delta}>> H^1(S,\mathcal I_{C/S}\otimes_S N_{S/X})
    @>>> 0,
  \end{CD}
\end{equation}
which is deduced from \eqref{ses:normal1}.
Then since $p_1$ is not surjective,
there exists a first order deformation $\tilde C$ of $C$ in $X$
not contained in any first order deformation $\tilde S$ of $S$ in $X$.
We need the following lemma for our proof of Theorem~\ref{thm:main1}.

\begin{lem}
  \label{lem:maximality}
  Suppose that $H^1(X,N_{(C,S)/X})=0$
  and $H^1(S,\mathcal I_{C/S}\otimes_S N_{S/X})\ne 0$.
  If the primary obstruction $\ob(\alpha)$ is nonzero in $H^1(C,N_{C/X})$
  for every global section $\alpha \in H^0(C,N_{C/X})\setminus \im p_1$,
  then $pr_1'$ is dominant in a neighborhood of $[C]$.
  If moreover $H^0(S,\mathcal I_{C/S}\otimes_S N_{S/X})=0$,
  then
  $$
  \dim_{[C]} \Hilb_P X = \dim_{(C,S)} \HF_{(P,Q)} X.
  $$
\end{lem}

\Proof
The proof is essentially same as that of \cite[Lemma~4.11]{Nasu4},
where $X$ and $S$ are assumed to be a smooth del Pezzo $3$-fold and
its smooth hyperplane section, respectively.
However we repeat the proof here for the reader's convenience.
We show that every small global deformation of $C$ in $X$ 
is contained in that of $S$ in $X$.
Let $T$ be a small neighborhood of $[C]$ in $\Hilb_P X$
and $T \hookrightarrow \Hilb_{P} X$ the embedding of $T$.
Then by base change, there exists a family $C_T \subset X \times T$
of curves in $X$ with a point $0 \in T$ such that $C_0 = C$.
Let $\Spec k[t]/(t^2) \rightarrow T$ be an
element of the Zariski tangent space of $T$ at $0$.
Then there exist a first order deformation
$\tilde C \rightarrow \Spec k[t]/(t^2)$ of $C$
and a global section $\alpha$ of $N_{C/X}$, correspondingly.
Then by assumption, $\alpha$ is contained in $\im p_1$,
and hence there exists a first order deformation 
$(\tilde C,\tilde S)$ of $(C,S)$ with $\tilde S \supset \tilde C$.
Since $\HF X$ is nonsingular at $(C,S)$,
%Then by infinitesimal lifting property of smoothness 
%(cf.~\cite[\S4]{Hartshorne10}), 
there exists a global deformation $(C_T,S_T)$ of $(C,S)$ 
over $T$ as a lift of $(\tilde C,\tilde S)$.
Thus $pr_1'$ is dominant near $[C]$.
If moreover $H^0(S,\mathcal I_{C/S}\otimes_S N_{S/X})=0$,
then $pr_1'$ is locally an embedding in a neighborhood of $(C,S)$
(cf.~\cite[\S2.2]{Nasu6}).
Thus we have
$$
\dim_{(C,S)} \HF_{P,Q} X=\dim \mathcal W_{C,S}=\dim pr_1'(\mathcal W_{C,S})
=\dim_{[C]} \Hilb_P X.\qed
$$ 

\begin{cor}
  \label{cor:maximality}
  Let $X=\mathbb P^3$,
  $S$ a smooth cubic surface in $\mathbb P^3$,
  and $C \subset S$ a smooth curve of degree $d>9$ and genus $g$.
  If $\ob(\alpha)\ne 0$ in $H^1(C,N_{C/\mathbb P^3})$
  for all $\alpha \not\in \im p_1$,
  then $\dim_{[C]} H(d,g)^{sc} = d+g+18$.
\end{cor}

\Proof
We see that $\mathcal I_{C/S}\otimes_S N_{S/\mathbb P^3}
\simeq -C-3K_S$ in $\Pic S$.
Since $-K_S.(-C-3K_S)=-d+9<0$,
we have $H^0(S,\mathcal I_{C/S}\otimes_S N_{S/\mathbb P^3})=0$.
Then we have proved the corollary by Lemma~\ref{lem:expected dimension}.
\qed

\subsection{Obstructedness criterion}
\label{subsec:obstructed}

In this section, we recall a result in \cite{Nasu8}
concerning primary obstructions to deforming curves on a $3$-fold.
Our proofs of Theorems~\ref{thm:main1} and \ref{thm:main2}
heavily depend on Theorem~\ref{thm:obstruction}.
We refer to \cite{Mukai-Nasu,Nasu4,Nasu5,Nasu8}
for more information about exterior components, 
infinitesimal deformations with pole, 
and also the proof of Theorem~\ref{thm:obstruction}.

Let $X$ be a projective $3$-fold and $C$ an irreducible curve on $X$.
We assume that there exists an intermediate surface $S$ such
that $C \hookrightarrow S \hookrightarrow X$ are regular embeddings.
Let $\alpha$ be a global section of $N_{C/X}$.
We consider a natural projection
$\pi_{C/S}: N_{C/X} \rightarrow N_{S/X}\big{\vert}_C$,
which induces maps
$H^i(C,N_{C/X}) \rightarrow H^i(C,N_{S/X}\big{\vert}_C)$ 
($i=0,1$) on their cohomology groups.
The images of $\alpha$ and $\ob(\alpha)$ 
in $H^i(C,N_{S/X}\big{\vert}_C)$ ($i=0,1$) by the induced maps
are called the {\em exterior component} of $\alpha$ and $\ob(\alpha)$
and denoted by $\pi_{C/S}(\alpha)$ and $\ob_S(\alpha)$, respectively.
By definition, if $\ob_S(\alpha)$ is nonzero then so is $\ob(\alpha)$.

We recall the definition of infinitesimal deformations with poles,
which was introduced in \cite{Mukai-Nasu}.
We are interested in a global section $\gamma$ of $N_{S/X}\big{\vert}_C$
such that 
$\gamma$ does not lift to a global section of $N_{S/X}$
but lifts to that of $N_{S/X}(E)$ ($:=N_{S/X}\otimes_S \mathcal O_S(E)$)
after admitting a pole along a divisor $E\ge 0$ on $S$.

\begin{dfn}
\label{dfn:infinitesimal deformation with pole}
Let $E$ be a nonzero effective Cartier divisor on $S$.
Then a rational section 
$\beta \in H^0(S,N_{S/X}(E))\setminus H^0(S,N_{S/X})$
is called an {\em infinitesimal deformation with pole}.
\end{dfn}

Here and later, 
for a sheaf $\mathcal F$ and a Cartier divisor $E$ on $S$,
we denote the sheaf $\mathcal F\otimes_S \mathcal O_S(E)$ by $\mathcal F(E)$.
When $\mathcal F$ is invertible,
we abuse notations and denote the invertible sheaf $\mathcal F(E)$ 
by $\mathcal F+E$.
If $E$ is effective then for every integer $i\ge 0$ there exists a natural map
\begin{equation}
  \label{map:admitting pole}
  H^i(S,\mathcal F) \rightarrow H^i(S,\mathcal F\otimes_S \mathcal O_S(E)).
\end{equation}
Given a cohomology class $\mathbf c$ in $H^i(S,\mathcal F)$,
we denote by $r(\mathbf c,E)$ the image of $\mathbf c$ by this map 
(and similarly for $\mathbf c\big{\vert}_C$ in $H^i(C,\mathcal F\big{\vert}_C)$).
Let $\mathbf k_C$ denote the extension class
in $\Ext^1_S(\mathcal O_C,\mathcal O_S(-C))$ of
the short exact sequence
\begin{equation}
\label{ses:k_C}
\begin{CD}
  0 @>>> \mathcal O_S(-C) @>>> \mathcal O_S
  @>>> \mathcal O_C @>>> 0
\end{CD}
\end{equation}
on $S$. When $E\ge 0$ is prime
and $\beta$ is a global section of $\mathcal F(E)$,
we call the restriction $\beta\big{\vert}_E$ of $\beta$ to $E$,
that is a global section of $\mathcal F(E)\big{\vert}_E$,
the {\em principal part} of $\beta$ {\em along} $E$.
The following lemma is a generalization of \cite[Lemma~3.1]{Nasu5}.

\begin{lem}
  \label{lem:with pole}
  Let $L$ be an invertible sheaf on $S$,
  $E$ a nonzero effective divisor on $S$
  not containing $C$ as its component,
  %$C\cap E$ are finitely many points,
  $\gamma$ a global section of $L\big{\vert}_C:=L\otimes_S \mathcal O_C$.
  Then
  \begin{enumerate}
    \item $r(\gamma,E)$ lifts to a global section $\beta$ of $L+E$ on $S$
    if and only if $r(\gamma,E) \cup \mathbf k_C=0$
    in $H^1(S,L+E-C)$.
    \item If $H^1(S,L+E-C)=0$ and $\gamma\cup \mathbf k_C\ne 0$,
    then there exist a triplet $(E',E_0,\beta)$ of 
    a subdivisor $E'\subset E$ on $S$,
    a prime divisor $E_0 \subset E'$,
    and a lift $\beta$ of $r(\gamma,E')$ in $H^0(S,L+E')$
    such that the principal part $\beta\big{\vert}_{E_0}$ 
    of $\beta$ along $E_0$ is nonzero and contained in the subgroup 
    $$
    H^0(E_0,(L+E'-C)\big{\vert}_{E_0})
    \subset H^0(E_0,(L+E')\big{\vert}_{E_0}).
    $$
  \end{enumerate}
\end{lem}

\Proof (1) follows from \cite[Lemma~3.1]{Nasu5}.
(Consider the first coboundary map 
of $\eqref{ses:k_C}\otimes_S L+E$, which is a map taking 
a cup product with $\mathbf k_C$.)
Since $E$ is nonzero and effective, 
there exist positive integers $k$, $m_i$ ($1 \le i \le k$) and 
prime divisors $E_i$ ($1 \le i \le k$) such that
$$
E=\sum_{i=1}^k m_i E_i.
$$
Since the reduction map $r(*,E)$ (cf.~\eqref{map:admitting pole})
and the cup product map $\cup \mathbf k_C$ are compatible,
we see that 
$$
r(\gamma,E')\cup \mathbf k_C = r(\gamma \cup \mathbf k_C,E')
$$
in $H^1(S,L+E'-C)$ for any subdivisor $E'\subset E$.
Thus by admitting to $\gamma$ a new pole along some $E_i$,
or increasing the order of poles along $E_i$,
we obtain a divisor $E' \subset E$ such that
\begin{enumerate}
  \renewcommand{\theenumi}{\alph{enumi}}
  \renewcommand{\labelenumi}{{\rm (\theenumi)}}
  \item $r(\gamma,E') \cup  \mathbf k_C= 0$ in $H^1(S,L+E'-C)$ and
  \label{item:liftable}
  \item $r(\gamma,E'-E_i) \cup  \mathbf k_C\ne 0$ in $H^1(S,L+E'-E_i-C)$.
  \label{item:not liftable}
\end{enumerate}
for some $1 \le i\le k$.
Then $r(\gamma,E')$ lifts to a global section $\beta$ of $L+E'$
by \eqref{item:liftable} but does not lift to that of $L+E'-E_i$
by \eqref{item:not liftable}.
Then by virtue of \cite[Lemma~3.1]{Nasu5},
$\beta\big{\vert}_{E_i}$ is nonzero in $H^0(E_i,(L+E')\big{\vert}_{E_i})$ 
and contained in the subgroup $H^0(E_i,(L+E'-C)\big{\vert}_{E_i})$.
\qed

\medskip

We recall a sufficient condition for $\ob_S(\alpha)$ to be nonzero.
Let $E_i$ ($1 \le i \le k$) be nonzero effective prime divisors on $S$
such that
\begin{enumerate}
  \item $E_i$ are mutually disjoint, 
  i.e.,~$E_i \cap E_j=\emptyset$ if $i \ne j$, and
  \item if $D$ and $D'$ are two effective divisors on $S$
  whose supports are contained in $\bigcup_{i=1}^k E_i$ and if 
  $D \le D'$, then the natural map
  \begin{equation}
    \label{map:divisors ass lines}
    H^1(S,D) \longrightarrow H^1(S,D')
  \end{equation}
  %inclusion $\mathcal O_S(D) \hookrightarrow \mathcal O_S(D')$
  is injective.
\end{enumerate}  

\begin{ex}
  \label{ex:lines and cubic}
  If $S$ is a del Pezzo surface and $E_i$ ($1 \le i \le k$) are
  mutually disjoint lines on $S$,
  then the map \eqref{map:divisors ass lines} is injective.
  In fact, since $D$ and $D'$ have supports on $\bigcup_{i=1}^k E_i$, 
  so does $E:=D'-D$.
  If $D\lneqq D'$ then $E$ is nonzero and effective.
  Then since $E_i$ are $(-1)$-curves on $S$,
  we see that $H^0(E,\mathcal O_E(D'))=0$,
  where $\mathcal O_E(D')\simeq \mathcal O_S(D')\otimes_S \mathcal O_E$.
  Thus the injectivity of \eqref{map:divisors ass lines}
  follows from the exact sequence
  $0 \rightarrow \mathcal O_S(D) \rightarrow \mathcal O_S(D')
  \rightarrow \mathcal O_E(D') \rightarrow 0$.
\end{ex}

Let $\gamma=\pi_{C/S}(\alpha)$ be the exterior component of $\alpha$.
We consider a divisor $E=\sum_{i=1}^k m_i E_i$ on $S$
with positive coefficients $m_i \in \mathbb Z_{>0}$ 
and assume that $C\ne E_i$ for any $i=1,\dots,k$ 
(then $C \cap E_i$ are finitely many points).
We assume furthermore that $r(\gamma,E)$ lifts to a section 
$\beta \in H^0(S,N_{S/X}(E))\setminus H^0(S,N_{S/X})$
(an infinitesimal deformation with pole), 
i.e.,~we have
\begin{equation}
  \label{eqn:lifting}
  r(\pi_{C/S}(\alpha),E)= \beta\big{\vert}_C
  \qquad
  \mbox{in}
  \qquad
  H^0(C,N_{S/X}(E)\big{\vert}_C).
\end{equation}
Let $\beta_i:=\beta\big{\vert}_{E_i}$ 
be the principal part of $\beta$ along $E_i$.
Then by assumption, $\beta_i$ is a global section of 
the invertible sheaf
$N_{S/X}(m_iE_i)\big{\vert}_{E_i}$ 
($\simeq N_{S/X}(E)\big{\vert}_{E_i}$) on $E_i$.
Moreover, by Lemma~\ref{lem:with pole}, 
$\beta_i$ is contained in the subgroup
$$
H^0(E_i,N_{S/X}(m_iE_i-C)\big{\vert}_{E_i})
\subset H^0(E_i,N_{S/X}(m_iE_i)\big{\vert}_{E_i}).$$
We illustrate the relations among $\alpha$, $\beta$ and $\beta_i$
in Figure~\ref{fig:relation}.
\begin{figure}[h]
$$
\begin{array}{ccccccc}
  H^0(C,N_{C/X}) & \ni & \alpha & & & & H^0(E_i,N_{E_i/X}(E)) \\
  \mapdown{\pi_{C/S}} && \big\downarrow && && \mapdown{\pi_{E_i/S}(E)} \\
  H^0(C,N_{S/X}\big{\vert}_C) & \ni & \gamma
  &  & H^0(S,N_{S/X}(E)) & \longrightarrow  
  & H^0(E_i,N_{S/X}(E)\big{\vert}_{E_i}) \\
 \mapdown{r} && \big\downarrow && \rotatebox{-90}{$\ni$} && \rotatebox{-90}{$\ni$}\\
 H^0(C,N_{S/X}(E)\big{\vert}_C) & \ni & r(\gamma,E) & \overset{res}\longmapsfrom & 
  \beta & \overset{res}\longmapsto&  \beta_i
\end{array}
$$
\caption{Relation among $\alpha$, $\beta$ and $\beta_i$}
\label{fig:relation}
\end{figure}

Let $\partial_{E_i}$ denote the coboundary map of the short exact sequence
\begin{equation}
  \label{ses:normal bundle of E} 
  \begin{CD}
    [0 @>>> N_{E_i/S} @>>> N_{E_i/X} @>{\pi_{E_i/S}}>> N_{S/X}\big{\vert}_{E_i}
      @>>> 0] \otimes_{E_i} \mathcal O_{E_i}(E)
  \end{CD}
\end{equation}
on $E_i$. Then $\partial_{E_i}(\beta_i)$ defines an element of
$H^1(E_i,N_{E_i/S}(E))$ ($\simeq H^1(E_i, (m_i+1)E_i)$).
The following theorem is a refinement of \cite[Theorem~1.1]{Nasu5},
which enables us to deduce the nonzero of $\ob_S(\alpha)$
from that of the cup product of $\partial_{E_i}(\beta_i)$ with $\beta_i$.

\begin{thm}[{cf.~\cite[Theorem~1]{Nasu8}}]\label{thm:obstruction}
  Suppose that $H^1(S,N_{S/X})=0$.
  Then the exterior component 
  $\ob_S(\alpha)$ of $\ob(\alpha)$ is nonzero 
  in $H^1(C,N_{S/X}\big{\vert}_C)$ if
  we have the following:
  \begin{enumerate}
    \item Let $\Delta:=C+K_X\big{\vert}_S-2E$ in $\Pic S$
    and let $E_{\red}:=\sum_{i=1}^k E_i$,
    i.e.,~the reduced part of $E$.
    Then the restriction map 
    $$
    H^0(S,\Delta) \overset{\vert_{E_{\red}}}\longrightarrow 
    H^0(E_{\red},\Delta\big{\vert}_{E_{\red}})
    $$
    to $E_{\red}$ is surjective, and 
    \label{item:surjective}
    \item There exists an integer $1 \le i \le k$ such that
    $\partial_{E_i}(\beta_i) \cup \beta_i \ne 0$,
    where the cup product is taken by the map
    $$
    H^1(E_i,(m_i+1)E_i)\times
    H^0(E_i,N_{S/X}(m_iE_i-C)\big{\vert}_{E_i}) \overset{\cup}{\longrightarrow}
    H^1(E_i,N_{S/X}((2m_i+1)E_i-C)\big{\vert}_{E_i}).
    $$
    \label{item:nonzero}
  \end{enumerate}
\end{thm}
%which is not correct as it stand, because its proof contains an error.
%The correction of the error and the proof of Theorem~\ref{thm:obstruction}
%will be given somewhere.
%For the proof of Theorem~\ref{thm:obstruction}, we refer to
%the postprint of the same paper in the arXiv (cf.~\cite{Nasu5}).
In the rest of this section, 
we assume that $X=\mathbb P^3$, 
$S \subset X$ is a smooth cubic surface,
and $E_i$ are lines on $S$.
The following lemma gives a sufficient condition 
for the cup product $\partial_{E_i}(\beta_i) \cup \beta_i$
considered in Theorem~\ref{thm:obstruction} to be nonzero.

\begin{lem}
  \label{lem:nonzero of cup product}
  Let $Z_i:=C\cap E_i$ be the scheme-theoretic intersection 
  of $C$ with $E_i$. If
  \begin{enumerate}
    \renewcommand{\theenumi}{\roman{enumi}}
    \renewcommand{\labelenumi}{{\rm [\theenumi]}}
    \item \label{item:principal}
    $\beta_i\ne 0$, equivalently, $\beta$ is not contained in
    $H^0(S,N_{S/\mathbb P^3}(E-E_i))$,
    
    \item \label{item:trivial}
    $(C.E_i)=3-m_i$, and
    
    \item \label{item:coboundary}
    If $m_i=1$ then $Z_i$ is a general member
    of a linear system $\Lambda:=|\mathcal O_{E_i}(2)|$
    on $E_i \simeq \mathbb P^1$,
  \end{enumerate}
  then the cup product $\partial_{E_i}(\beta_i) \cup \beta_i$ 
  is nonzero.
\end{lem}
\Proof
Since $N_{S/\mathbb P^3}\simeq \mathcal O_S(3)$
and $E_i$ is a $(-1)$-curve, 
we see that $N_{S/\mathbb P^3}(m_iE_i)\vert_{E_i}$
is an invertible sheaf on $E_i \simeq \mathbb P^1$ of degree $3-m_i$.
Then it has a nonzero section by [\ref{item:principal}],
and thereby we obtain $3-m_i\ge 0$.
This implies that $m_i=1,2$ or $3$.
Then it follows from the condition [\ref{item:trivial}] that 
$N_{S/\mathbb P^3}(m_iE_i-C)\vert_{E_i}$ is a trivial sheaf.
Therefore, taking a cup product with $\beta_i$ is just 
a multiplication by a nonzero scalar.
Hence for the proof, it suffices to prove
that $\partial_{E_i}(\beta_i)\ne 0$.
We consider the map 
$$
\pi_{E_i/S}(E): H^0(E_i,N_{E_i/\mathbb P^3}(m_iE_i))
\longrightarrow H^0(E_i,N_{S/\mathbb P^3}(m_iE_i)\big{\vert}_{E_i}),
$$
which is induced by a sheaf homomorphism 
$\pi_{E_i/S} \otimes_{E_i} \mathcal O_{E_i}(E)$
in \eqref{ses:normal bundle of E}.
We see that this map is zero if $m_i>1$ and
not surjective if $m_i=1$,
because
$N_{E_i/\mathbb P^3} \simeq \mathcal O_{\mathbb P^1}(1)^{\oplus 2}$,
$N_{S/\mathbb P^3}\big{\vert}_{E_i}\simeq \mathcal O_{\mathbb P^1}(3)$ and
$\mathcal O_{E_i}(E_i)\simeq \mathcal O_{\mathbb P^1}(-1)$
on $E_i \simeq \mathbb P^1$.
Thus if $m_i >1$ then we are done. Suppose that $m_i=1$.
Then by the condition [\ref{item:coboundary}],
$Z_i$ is a finite subscheme of $E_i$ of length $2$.
Given an invertible sheaf $\mathcal L$ and its global section $\gamma$,
we denote by $\div_0(\gamma)$ the divisor of zero of $\gamma$.
Then by \cite[Lemma~3.1]{Nasu5},
$\beta_i$ is contained in the subgroup
$H^0(E_i,N_{S/\mathbb P^3}(E_i-C)\big{\vert}_{E_i})
\subset H^0(E_i,N_{S/\mathbb P^3}(E_i)\big{\vert}_{E_i})$.
Therefore, as a section of the sheaf
$N_{S/\mathbb P^3}(E_i)\big{\vert}_{E_i}\simeq \mathcal O_{\mathbb P^1}(2)$
on $E_i \simeq \mathbb P^1$, we have $\div_0 (\beta_i)=Z_i$.
If $\partial_{E_i}(\beta_i)=0$, then 
$Z_i$ is contained in the linear subsystem
$$
\left\{\div_0(\gamma)\bigm|  \gamma \in \im \pi_{E_i/S}(E_i) \right\}
\subsetneqq |N_{S/\mathbb P^3}(E_i)\big{\vert}_{E_i}|
$$
of codimension one, thereby contradicting the genericity of $Z_i$
as mentioned in the condition [\ref{item:coboundary}].
Thus we conclude that $\partial_{E_i}(\beta_i)\ne 0$.\qed

\medskip
The next lemma will be used to prove 
that the condition [\ref{item:coboundary}]
in Lemma~\ref{lem:nonzero of cup product} 
is satisfied in the case where $m_i=1$, i.e.,~$C.E_i=2$.

\begin{lem}
  \label{lem:general member}
  Let $E$ be a line and $D$ a nef divisor on $S$.
  If $D\not\sim m(-K_S-E)$ for any integer $m$,
  then we have $H^1(S,D-E)=0$, and in particular
  the rational map
  $|D| \dashrightarrow |\mathcal O_E(D)|$
  sending a curve $C \in |D|$ to $Z:=C\cap E$ is dominant.
\end{lem}
\Proof
We note that $q:=-K_S-E$ is the class of conics on $S$ residual to $E$.
Put $L:=D+q$ in $\Pic S$. Since $D$ is nef, so is $L$.
Since $L$ is not composed with pencils, $L$ is also big.
Thus $H^1(S,D-E)\simeq H^1(S,K_S+L)=0$ as a consequence of 
Kawamata-Viehweg vanishing theorem.
Then the lemma follows from the exact sequence
\begin{equation}
  \label{ses:k_E}
  \begin{CD}
    0 @>>> \mathcal O_S(-E) @>>> \mathcal O_S @>>> \mathcal O_E 
    @>>> 0]\otimes_S \mathcal O_S(D). \qed
\end{CD}
\end{equation}

\section{Proof of Theorems~\ref{thm:main1} and \ref{thm:main2}}
\label{sec:proof}

In this section, we prove Theorems~\ref{thm:main1} and \ref{thm:main2}.
Let $S \subset \mathbb P^3$ be a smooth cubic surface,
$C \subset \mathbb P^3$ a curve contained in $S$.
We define an invertible sheaf $L$ on $S$ by
$L:=\mathcal O_S(C)\otimes_S N_{S/\mathbb P^3}^{-1}$.
Then $L \sim C+3K_S$ in $\Pic S$.
We will see that the two cohomology groups
$H^1(S,-L)$ and $H^2(S,-L)$ on $S$ are important
for studying the deformations of $C$ in $\mathbb P^3$.

Since $H^1(C,N_{C/S})=0$, as we saw in \S\ref{sec:intro},
the cohomology group $H^1(C,N_{S/\mathbb P^3}\big{\vert}_C)$
contains every obstruction
to deforming $C$ in $\mathbb P^3$.
Since $H^i(S,N_{S/\mathbb P^3})=0$ for all $i>0$,
it follows from the exact sequence 
$\eqref{ses:k_C}\otimes_S N_{S/\mathbb P^3}$ that
\begin{equation}
  \label{isom:obstruction space}
  H^1(C,N_{S/\mathbb P^3}\big{\vert}_C)\simeq H^2(S,-L).
\end{equation}
The following lemma shows that
if $H^1(S,-L)\ne 0$ and $\chi(S,-L)\ge 0$, 
then the obstruction space
\eqref{isom:obstruction space} is nonzero
by $H^2(S,-L) \simeq H^0(S,L+K_S)^{\vee}$.
\begin{lem}
  \label{lem:effectivity of L+K_S}
  Suppose that $C$ is of degree $d>9$ and genus $g\ge 3d-18$.
  Then
  \begin{enumerate}
    \item $H^0(S,-L)=0$,
    \item $\chi(S,-L)\ge 0$ and
    \item If $H^1(S,-L)\ne 0$ then 
    $L+K_S\ge 0$, $L$ is big and not nef.
  \end{enumerate}
\end{lem}
\Proof
(1) follows from $L.K_S=9-d<0$.
By Riemann-Roch theorem on $S$, we have
$$
\chi(S,-L)=(C+3K_S)(C+4K_S)/2+1=g-3d+18
$$
and hence we obtain (2) by assumption.
This implies that if $H^1(S,-L)\ne 0$ then $H^2(S,-L)\ne 0$.
Then it follows from Lemma~\ref{lem:not nef} that
$L+K_S$ is effective and $L$ is not nef.
Finally $L$ is big by Lemma~\ref{lem:big}.
\qed

\medskip

We next relate $H^1(S,-L)$
to a tangent map on the Hilbert-flag scheme.
By Lemma~\ref{lem:fano-fano}, we note that
$H^i(\mathbb P^3,N_{(C,S)/\mathbb P^3})=0$ for all $i>0$.
This implies that
the Hilbert-flag scheme $\HF \mathbb P^3$ of $\mathbb P^3$ is nonsingular
at $(C,S)$ of expected dimension 
$\chi(\mathbb P^3,N_{(C,S)/\mathbb P^3})=d+g+18$
(cf.~Lemma~\ref{lem:expected dimension}).
Let $\Hilb^{sc} \mathbb P^3$ denote the Hilbert scheme of 
smooth connected curves in $\mathbb P^3$,
and let $\HF^{sc} \mathbb P^3 \subset \HF \mathbb P^3$ 
denote the subscheme parametrising pairs $(C',S')$ of a curve 
$C' \in \Hilb^{sc} \mathbb P^3$ and 
a surface $S'$ containing $C'$, i.e.,~we define by
$\HF^{sc} \mathbb P^3:=pr_1^{-1}(\Hilb^{sc} \mathbb P^3)$,
where $pr_1: \HF \mathbb P^3 
\rightarrow \Hilb \mathbb P^3$ is the first projection.
Then by \eqref{ses:first tangent map} the cokernel of the tangent map
\begin{equation}
  \label{map:tangent map}
  p_1: H^0(\mathbb P^3,N_{(C,S)/\mathbb P^3}) 
  \longrightarrow H^0(C,N_{C/\mathbb P^3})
\end{equation}
of $pr_1$ at $(C,S)$ is isomorphic to 
$H^1(S,N_{S/\mathbb P^3}(-C))=H^1(S,-L)$.
The following lemma immediately follows
from Corollary~\ref{cor:vanishing and base component} and 
Lemma~\ref{lem:abnormality}.

\begin{lem}
  \label{lem:quadratically normal}
  Suppose that $L\ge 0$ and $L$ is not nef.
  Then the fixed part $F$ of $|L|$ is given by
  $$
  F=m_1E_1 +\cdots + m_kE_k,
  $$
  where $E_i$ ($1 \le i \le k$) are mutually disjoint lines on $S$
  such that $L.E_i<0$ and $m_i=-L.E_i$. 
  Here we have $m_i\le 3$ for all $i$.
  If moreover $L^2>0$ then we have the followings:
  \begin{enumerate}
    \item $h^1(S,-L)=h^0(F,\mathcal O_F)$, and
    \item $C$ is quadratically normal (resp. linearly normal)
    if and only if $m_i=1$ (resp. $1 \le m_i\le 2$)
    for all $i$.
  \end{enumerate}
\end{lem}

\medskip

\paragraph{\bf Proof of Theorem~\ref{thm:main1}}

We assume that $C$ satisfies the hypothesis of Conjecture~\ref{conj:Kleppe}.
Then $H^1(S,-L)\simeq H^1(\mathbb P^3,\mathcal I_C(3))\ne 0$
by \eqref{isom:abnormality}.
Thereby the tangent map $p_1$ defined above is not surjective.
Let $\alpha$ be a global section of $N_{C/\mathbb P^3}$.
Then by lemma~\ref{lem:maximality} 
(or more directly by Corollary~\ref{cor:maximality}),
it suffices to prove that the primary obstruction
$\ob(\alpha)$ of $\alpha$ is nonzero 
if $\alpha$ is not contained in the image of $p_1$.
It follows from Lemma~\ref{lem:effectivity of L+K_S} that
$\chi(S,-L)\ge 0$, $L+K_S\ge 0$, $L$ is not nef and $L^2>0$.

Suppose now that $C$ is quadratically normal.
Then Lemma~\ref{lem:quadratically normal} shows that
the fixed part $F$ of $|L|$ is given by
$$
F=E_1 +\cdots + E_k,
$$
where $k=h^1(S,-L)$ and $E_i$ ($1 \le i \le k$) 
are lines on $S$ mutually disjoint.
Then $L-F$ is clearly nef, and also big by $(L-F)^2=L^2-F^2>L^2>0$.
Thus we see that $H^1(S,-L+F)=0$.
Let $\mathbf k_C$ be the extension class defined by \eqref{ses:k_C}
and $\gamma:=\pi_{C/S}(\alpha)$ the exterior component of $\alpha$
(see \S\ref{subsec:obstructed} for the definition).
We note that the map $\delta$ in \eqref{ses:first tangent map},
which is the first coboundary map of \eqref{ses:normal1},
factors through the coboundary map
$\cup \mathbf k_C$ of $\eqref{ses:k_C}\otimes_S N_{S/\mathbb P^3}$
(cf.~\cite[Lemma~2.2]{Nasu6}).
Then since $\alpha$ is not contained in $\im p_1$,
the cup product
$\gamma\cup \mathbf k_C$ is nonzero in $H^1(S,-L)$.
We note that $H^1(S,N_{S/\mathbb P^3}(F-C))\simeq H^1(S,-L+F)=0$.
Then by Lemma~\ref{lem:with pole},
admitting to $\gamma$ some poles along $F\cap C$,
the section $\gamma$ (in fact $r(\gamma,F)$) on $C$
lifts to a global section of $N_{S/\mathbb P^3}(F)$ on $S$,
i.e., an infinitesimal deformation with pole
(cf.~Definition~\ref{dfn:infinitesimal deformation with pole}).
More precisely, by the same lemma,
there exist %a subdivisor $E \subset F$ on $S$,
a prime divisor $E_i \subset F$ on $S$ ($1 \le i \le k$)
and a lift $\beta \in H^0(S,N_{S/\mathbb P^3}(F))$ of $r(\gamma,F)$
such that the principal part
$\beta_i:=\beta\big{\vert}_{E_i}$ of $\beta$ along $E_i$ is nonzero
in $H^0(E_i,N_{S/\mathbb P^3}(F)\big{\vert}_{E_i})\simeq
H^0(E_i,N_{S/\mathbb P^3}(E_i)\big{\vert}_{E_i})$.

Now we check that the two conditions
\eqref{item:surjective} and \eqref{item:nonzero}
of Theorem~\ref{thm:obstruction} are both satisfied.
Let us define a divisor $\Delta$ on $S$ as in the theorem.
Then $\Delta=C+K_{\mathbb P^3}\big{\vert}_S-2F\sim L+K_S-2F$.
The Serre duality shows that 
$H^1(S,\Delta-F) \simeq H^1(S,3F-L)^{\vee}$
and the last cohomology group is zero by Lemma~\ref{lem:crutial}.
Therefore the restriction map
$H^0(S,\Delta)\rightarrow H^0(F,\Delta\big{\vert}_{F})$ is 
surjective. Thus \eqref{item:surjective} is satisfied.

To check the condition \eqref{item:nonzero} of
Theorem~\ref{thm:obstruction},
we prove that the three conditions 
[\ref{item:principal}], [\ref{item:trivial}] 
and [\ref{item:coboundary}] of
Lemma~\ref{lem:nonzero of cup product} are all satisfied.
[\ref{item:principal}] is clear. 
[\ref{item:trivial}] follows from $m_i=1$ and $C.E_i=2$.
Since $C$ is a general member of the $3$-maximal family $W$,
so is $Z_i:=C\cap E_i$ in $|\mathcal O_{E_i}(2)|$
on $E_i \simeq \mathbb P^1$ by Lemma~\ref{lem:general member}.
Thus [\ref{item:coboundary}] follows.
Then the cup product
$\partial_{E_i}(\beta_i)\cup \beta_i$ 
considered in Theorem~\ref{thm:obstruction}~\eqref{item:nonzero} 
is nonzero. Thereby we have proved Theorem~\ref{thm:main1}.
\qed

\begin{rmk}
  \label{rmk:known ranges}
  In this remark, we collect some known results related to
  Conjecture~\ref{conj:Kleppe}.
  Kleppe~\cite{Kleppe87} proved the conjecture is true 
  in the range of the $(d,g)$-plane:
  $g>-1+(d^2-4)/8$ for $14\le d\le 17$ and $g>7+(d-2)^2/8$ for $d\ge 18$.
  Later, Ellia~\cite{Ellia87} proved the conjecture
  in the wider range: 
  $g>G(d,5)$ for $d\geq 21$, 
  where $G(d,5)$ denotes the maximal genus of curves of degree $d$
  not contained in a quartic surface and $G(d,5)\approx d^2/10$
  for $d\gg 0$ (cf.~\cite{Gruson-Peskine82}).
  It has been proved in \cite{Nasu1} that
  the conjecture is true if $h^1(\mathbb P^3,\mathcal I_C(3))=1$
  (that is the case $b_6=2$ and $b_5\ge 3$) by a method of this paper.
  Recently in the appendix of \cite{Kleppe-Ottem15}
  and more recently in \cite{Kleppe17},
  Kleppe has further extended the known range of $(d,g)$ where
  Conjecture~\ref{conj:Kleppe} holds to be true
  by a method of \cite{Kleppe87} together with a result in \cite{Ellia87},
  but his result does not cover our result.
  %and has made a progress in proving the conjecture.
  It is notable that his result shows that the conjecture is true
  for some classes of quadratic non-normal curves $C$
  (with $b_6=1$, $b_5\ge 5$ and satisfying some further assumptions)
  (cf.~\cite[Theorem~A.3]{Kleppe-Ottem15}).
  As far as we know,
  every proof that has been known so far is partial,
  and Conjecture~\ref{conj:Kleppe} is still open
  (in the case where $C \subset \mathbb P^3$ 
    is quadratically non-normal).
\end{rmk}

In the rest of this section, we prove Theorem~\ref{thm:main2}.

\medskip

\paragraph{\bf Proof of Theorem~\ref{thm:main2}}

Let $C$ satisfy the assumption of the theorem.
For the proof, it suffices to show that
there exists a global section $\alpha$ of $N_{C/\mathbb P^3}$
such that $\ob(\alpha)$ (or its exterior component $\ob_S(\alpha)$) 
is nonzero. 

Let $m:=-L.E$ and suppose that $1 \le m \le 3$.
If $1 \le j \le m$ then $(-L+jE).E=m-j\ge 0$. Therefore,
there exists a sequence 
$$
H^1(S,-L) \twoheadrightarrow H^1(S,-L+E) \twoheadrightarrow \dots 
\twoheadrightarrow H^1(S,-L+mE)
$$
of natural surjective maps. 
Since $L+K_S$ is effective by assumption,
so are $L$ and $L-mE$ by Corollary~\ref{cor:vanishing and base component}.
Since $L\not\sim mE$, we have $H^0(S,-L+mE)=0$.
Then since the invertible sheaf $\mathcal O_E(-L+mE)$ 
on $E\simeq \mathbb P^1$ is trivial, 
we deduce from the exact sequence 
\eqref{ses:k_E} (for $D=-L+mE$) that 
$$
h^1(S,-L+(m-1)E)-h^1(S,-L+mE)=1.
$$
Therefore there exists an element $\xi$ of $H^1(S,-L)$
such that $r(\xi,(m-1)E)\ne 0$ in $H^1(S,-L+(m-1)E)$
and $r(\xi,mE)=0$ in $H^1(S,-L+mE)$.
Here and later,
we use the same notation $r(*,D)$ in \S\ref{subsec:obstructed}
for a divisor $D \ge 0$ on $S$.
It follows from the exact sequence \eqref{ses:first tangent map}
that there exists a global section $\alpha$ of 
$N_{C/\mathbb P^3}$ such that $\delta(\alpha)=\xi$,
where $\delta$ is the first coboundary map of \eqref{ses:normal1}.
Let $\gamma:=\pi_{C/S}(\alpha)$ denote 
the exterior component of $\alpha$ (see \S\ref{subsec:obstructed})
and let $\mathbf k_C$ be the extension class of \eqref{ses:k_C}.
Then $\gamma \cup \mathbf k_C=\xi$
as in the proof of Theorem~\ref{thm:main1}.
Moreover by the choice of $\xi$ and Lemma~\ref{lem:with pole},
there exists a global section $\beta$ of 
$N_{S/\mathbb P^3}(mE)$ such that 
$\beta\big{\vert}_C=r(\gamma,mE)$
and the principal part $\beta\big{\vert}_E$ of $\beta$ along $E$
defines a nonzero global section of
$N_{S/\mathbb P^3}(mE-C)\big{\vert}_E$ on $E$.

Let us define a divisor $\Delta$ on $S$
by $\Delta:=C+K_{\mathbb P^3}\big{\vert}_S-2mE$ 
as in Theorem~\ref{thm:obstruction}.
We check that $\Delta$ and
$\beta$ (or $\beta\big{\vert}_E$)
satisfy the two assumptions
\eqref{item:surjective} and \eqref{item:nonzero} of the theorem.
Let $\varrho$ denote the restriction map defined by \eqref{map:restriction}.
If $m=1$ then $\Delta=L+K_S-2E$ is effective
by Corollary~\ref{cor:vanishing and base component}. Thus we see that
$\varrho$ is surjective for $m=1$ by $\Delta.E=0$
and Lemma~\ref{lem:surjective} below, 
and also for $2 \le m \le 3$ by assumption.
Thus \eqref{item:surjective} is satisfied.
To prove that the cup product
$\partial_{E}(\beta\big{\vert}_E) \cup \beta\big{\vert}_E$
considered in Theorem~\ref{thm:obstruction} is nonzero, 
we again apply Lemma~\ref{lem:nonzero of cup product}.
We have already seen that $\beta\big{\vert}_E\ne 0$ 
(cf.~[\ref{item:principal}]).
It is also clear that $m(:=-L.E)=3-C.E$ (cf.~[\ref{item:trivial}]).
Finally, if $m=1$, i.e.,~$C.E=2$,
replacing $C$ with a general member $C'$ of $|C|$, we can assume that
the intersection $Z=C\cap E$, that is a divisor on $E\simeq \mathbb P^1$ 
of degree $2$, is general in $|\mathcal O_E(2)|$
by Lemma~\ref{lem:general member} (cf.~[\ref{item:coboundary}]).
In fact, if $C'$ is obstructed in $\mathbb P^3$,
then so is $C$ by upper semicontinuity.
Thereby we have obtained all the desired properties 
of $\beta\big{\vert}_E$ enough for proving that
its cup product with $\partial_{E}(\beta\big{\vert}_E)$
is nonzero. Then by Theorem~\ref{thm:obstruction},
we have completed the proof of Theorem~\ref{thm:main2}.
\qed

\begin{lem}
  \label{lem:surjective}
  Let $E$ be a line on $S$ 
  and $\Delta$ a divisor on $S$ such that $n:=\Delta.E\ge 0$.
  If there exists a conic $q$ on $S$
  such that $q.E=1$ and $\Delta-nq\ge 0$,
  then the restriction map $\varrho$ in \eqref{map:restriction}
  is surjective.
\end{lem}
\Proof
Let $q':=-K_S-E$. Then by $q'.E=2$,
we have $nq \not\sim mq'$ for any integer $m$.
Then it follows from Lemma~\ref{lem:general member} that
$H^0(S,nq) \rightarrow H^0(E,nq\big{\vert}_E)$ is surjective.
Since $\Delta-nq$ is effective,
$|\Delta|$ contains $|nq|$ as a linear subsystem.
We note that $\mathcal O_E(nq)\simeq \mathcal O_E(\Delta)$ by degree.
Since the restriction of $\varrho$ to $H^0(S,nq)$ is surjective, 
so is $\varrho$.
\qed

\begin{rmk}
  \label{rmk:obstructed}
  Some special cases of Theorem~\ref{thm:main2} were 
  also proved in \cite{Dolcetti-Pareschi88} ($m=3$) and \cite{Nasu1} ($m=1$).
  The same conclusion was proved in \cite{Nasu1} 
  under the assumption that $F=\Bs|L|$ 
  is a (single) line (cf.~\cite[Proposition~3.1]{Nasu1}).
  Dolcetti and Pareschi~\cite{Dolcetti-Pareschi88}
  proved that if $d \ge 21$ and $G(d,5) < g\le d^2/8-d/2+1$, 
  then every linearly non-normal
  curve $C \in H(d,g)^{sc}$ lying on a smooth cubic surface
  belongs to a non-reduced component of $H(d,g)^{sc}$
  of dimension $d+g+20$
  (hence $C$ is obstructed in $\mathbb P^3$),
  whose general members are linearly non-normal curves
  lying on a quartic surface with a double conic
  (cf.~\cite[Theorem~2.1]{Dolcetti-Pareschi88}).
  Such curves are generic projections of 
  curves lying a smooth quartic del Pezzo surface in $\mathbb P^4$,
  and this fact was first pointed out by Ellia~\cite{Ellia87}.
  See \cite{Kleppe-Ottem15} for examples of obstructed curves
  with $m=2$ (cf.~Remark~\ref{rmk:known ranges}).
\end{rmk}
  
\section{Examples}
\label{sec:examples}

In this section, we consider applications of
Theorems~\ref{thm:main1} and \ref{thm:main2}.
We first look at applications of Theorem~\ref{thm:main1}.
We give two series of $3$-maximal families of space curves
satisfying the assumption of Conjecture~\ref{conj:Kleppe}.

Let $(a;b_1,\dots,b_6)$ be a $7$-tuple of integers satisfying 
the set of conditions \eqref{conds:standard-prescribed-nef-big},
and let $W:=W(a;b_1,\dots,b_6)$ be
the irreducible closed subset of $H(d,g)^{sc}$ associated to it
(cf.~Definition~\ref{dfn:3-maximal}).
We denote  by $C$ a general member of $W$.
Then if $d>9$, $W$ becomes a $3$-maximal family of $H(d,g)^{sc}$
(cf.~\S\ref{sec:intro})
and $C$ is contained in a unique smooth cubic surface $S$.
In Examples~\ref{ex:non-reduced1} and \ref{ex:non-reduced2} below, 
we have $b_6=2$ and 
$C$ is quadratically normal by Lemma~\ref{lem:abnormality}.
Then by virtue of Theorem~\ref{thm:main1},
$W$ becomes a component of $(H(d,g)^{sc})_{\red}$ of $d+g+18$,
and moreover $H(d,g)^{sc}$ is generically non-reduced along $W$.
Thus in these example, the Hilbert scheme
$H(d,g)^{sc}$ is highly singular along $W$.
In fact, at the generic point $C \in W$,
the tangential dimension of the Hilbert scheme is greater
than its dimension as a scheme
by $h^1(\mathbb P^3,\mathcal I_C(3))$
(cf.~\eqref{eqn:tangential codimension}).
We compute the two numbers $h^1(\mathbb P^3,\mathcal I_C(3))$
and $h^1(C,\mathcal O_C(3))$, where the latter
represents the dimension of the obstruction space
$H^1(C,N_{C/\mathbb P^3})$ of $H(d,g)^{sc}$ at $[C]$.
It can be computed by the formula
\begin{equation}
  \label{eqn:obstruction space}
  h^1(C,N_{C/\mathbb P^3})=h^1(C,\mathcal O_C(3))=
  h^0(S,C+4K_S),
\end{equation}
which is deduced from \eqref{isom:obstruction space}
and the Serre duality. In the following examples, 
$F$ denotes the fixed part of the linear system $|C+3K_S|$ on $S$.

\begin{ex}
  \label{ex:non-reduced1}
  Let $\lambda$ be a non-negative integer and let
  $$
  W=W(\lambda+14;2,2,2,2,2,2) \subset H(d,g)^{sc}.
  $$
  Then we have $d=3(\lambda+10)$ and $g=(\lambda+16)(\lambda+9)/2$,
  thus $\dim W=d+g+18=(\lambda+16)(\lambda+15)/2$.
  By Lemma~\ref{lem:vanishing and base component}, we have
  $F= (0;-1,-1,-1,-1,-1,-1)=\sum_{i=1}^6 \mathbf e_i$.
  Therefore, by the method of computations used in Example \ref{ex:fixed part},
  we see that
  $h^1(\mathbb P^3,\mathcal I_C(3))=h^0(F,\mathcal O_F)=6$.
  It follows from \eqref{eqn:obstruction space} that
  $h^1(C,\mathcal O_C(3))=h^0(S,C+4K_S-2F)
  =h^0(\mathbb P^2,\mathcal O_{\mathbb P^2}(\lambda+2))
  =(\lambda+4)(\lambda+3)/2$.
\end{ex}

\begin{ex}
  \label{ex:non-reduced2}
  Let $\lambda\ge 0$ and let 
  $$
  W=W(\lambda+17;\lambda+8,7,2,2,2,2) \subset H(d,g)^{sc}.
  $$
  Then $d=2(\lambda+14)$ and $g=8\lambda+67$, thus $\dim W=10\lambda+113$.
  We note that $F$ consists of $5$ disjoint lines on $S$ by
  $F= (1;1,1,-1,-1,-1,-1)=(\mathbf l-\mathbf e_1-\mathbf e_2)
  +\sum_{i=3}^6 \mathbf e_i$,
  where $\mathbf l$ is the class of the pullback of lines in $\mathbb P^2$
  (cf.~Example \ref{ex:fixed part}).
  This implies that $h^1(\mathbb P^3,\mathcal I_C(3))=5$.
  Moreover, by using the formula \eqref{eqn:obstruction space} again,
  we compute that $h^1(C,\mathcal O_C(3))=h^0(S,C+4K_S-2F)
  =h^0(S,(\lambda+3)\mathbf l - (\lambda+2)\mathbf e_1 -\mathbf e_2)
  %(\lambda+5)(\lambda+4)/2-(\lambda+3)(\lambda+2)/2-1
  =2\lambda+6$.
\end{ex}

We next look at applications of Theorem~\ref{thm:main2}.
Let $S$ be a smooth cubic surface in $\mathbb P^3$.
We fix a line $E$ on $S$
and denote by $\varepsilon$ the blow-down $S\rightarrow S'$ of $E$.

\begin{prop}
  \label{prop:obstructed}
  Let $0 \le k \le 2$ be any integer and
  $q$ a conic on $S$ such that $q.E=1$
  and $D'$ a nef divisor on $S'$.
  Let $D$ be a divisor on $S$ defined by
  $$
  D=-4K_S+2(3-k)E+(2-k)q+\varepsilon^*D',
  $$
  and let $\Lambda:=|D|$ be the linear system on $S$ spanned by $D$. Then
  \begin{enumerate}
    \item \label{item:smooth and connected}
    every general member $C$ of $\Lambda$ is smooth and connected,
    \item \label{item:intersection number with line}
    $C.E=k$ and
    \item \label{item:obstructed}
    $C$ is obstructed in $\mathbb P^3$.
  \end{enumerate}
\end{prop}
\Proof
Since $D.E=4-2(3-k)+2-k=k\ge 0$, $D$ is nef and hence $\Lambda$ 
is base point free (cf.~\S\ref{subsec:del Pezzo}). 
Then \eqref{item:smooth and connected} 
follows from Bertini's theorem,
\eqref{item:intersection number with line} from $D.E=k$.
We put $L:=C+3K_S$. Then $m:=-L.E=3-k$ and we have $1\le m \le 3$.
Let $\Delta:=L+K_S-2mE=C+4K_S-2mE=(2-k)q+\varepsilon^* D'$.
Then $\Delta.E=2-k\ge 0$.
Since $\Delta-(\Delta.E)q=\varepsilon^*D' \ge 0$, 
it follows from Lemma~\ref{lem:surjective}
that the restriction map $\varrho$ in \eqref{map:restriction} 
is surjective (for all $0 \le k\le 2$).
Thus \eqref{item:obstructed} follows from Theorem~\ref{thm:main2}.
\qed

\medskip

By taking $E$ and $q$ in Proposition~\ref{prop:obstructed}
as $E=\mathbf e_6$ and $q=\mathbf l-\mathbf e_6$, respectively,
we obtain the following example.

\begin{ex}
  \label{ex:obstructed}
  Let $0 \le k \le 2$ be an integer 
  and let $(a;b_1,\dots,b_5)$ be a $6$-tuple of integers
  satisfying $a\ge b_1+b_2+b_3$ and $b_1\ge b_2\ge \dots \ge b_5 \ge 0$.
  Since the invertible sheaf
  $\mathcal O_{S'}(a;b_1,\dots,b_5)$ on $S'$ is nef 
  (cf.~Lemma~\ref{lem:lines and nefness}),
  it follows from Proposition~\ref{prop:obstructed} that
  every general member $C$ of the linear system
  $$
  |\mathcal O_S(14-k+a;b_1+4,b_2+4,b_3+4,b_4+4,b_5+4,k)|
  $$
  on $S$ is a smooth connected curve in $\mathbb P^3$
  with $C.\mathbf e_6=k$. 
  Moreover, $C$ is obstructed in $\mathbb P^3$.
\end{ex}

Theorem~\ref{thm:main2} can be applied to
determinations of the dimension of the Hilbert scheme $H(d,g)^{sc}$.
We first recall a result due to Kleppe.
\begin{thm}[{\cite[Theorem~1.1]{Kleppe87}}]
  \label{thm:codimension}
  Let $W$ be a $3$-maximal family in $H(d,g)^{sc}$
  whose general member $C$ is contained in a smooth cubic surface.
  If $d>9$ (or $h^0(\mathbb P^3, \mathcal I_C(3))=1$), then we have
  $$
  h^1(\mathbb P^3,\mathcal I_C(3))-h^1(C,\mathcal O_C(3))
  \le
  \dim_{[C]} H(d,g)^{sc}-\dim W
  \le 
  h^1(\mathbb P^3,\mathcal I_C(3)),
  $$
  where the inequality to the right is strict
  if and only if $C$ is obstructed in $\mathbb P^3$.
\end{thm}
See \cite[Theorem~2.4]{Nasu6} for a generalization of this theorem.
Theorems~\ref{thm:codimension} and \ref{thm:main2}
allow us to determine the dimension of $H(d,g)^{sc}$ at $[C]$
in the case where $h^1(C,\mathcal O_C(3))=1$. 

\begin{prop}
  \label{prop:codimension}
  Let $C \subset \mathbb P^3$ be a smooth connected curve 
  of degree $d$ and genus $g$ lying on a smooth cubic surface $S$.
  \begin{enumerate}
    \item \label{item:dimension}
    If $h^0(\mathbb P^3, \mathcal I_C(3))=1$, $h^1(C,\mathcal O_C(3))=1$ and
    $C$ is obstructed, then
    \begin{equation}
      \label{eqn:dimension2}
      \dim_{[C]} H(d,g)^{sc} =d+g+17+h^1(\mathbb P^3,\mathcal I_C(3)).
    \end{equation}
    \item \label{item:dimension2}
    Suppose that $C$ is a member of the linear system
    $$
    |\mathcal O_S(12;b_1,b_2,\dots,b_6)|
    $$
    on $S$ with $b_i$ satisfying $0 \le b_i\le 4$ for all $i$.
    If $b_j=2$ for some $1 \le j \le 6$, then we have 
    \eqref{eqn:dimension2}.
  \end{enumerate}
\end{prop}
\Proof
\eqref{item:dimension} follows from Theorem~\ref{thm:codimension}
and $\dim W=d+g+18$.
We prove \eqref{item:dimension2}.
We note that $L+K_S=C+4K_S=\sum_{i=1}^6 (4-b_i)\mathbf e_i\ge 0$.
Since $-L.\mathbf e_j=1$, $C$ is obstructed by Theorem~\ref{thm:main2}.
Moreover, since $C+4K_S$ is a sum of lines on $S$, 
we see that 
$h^1(C,\mathcal O_C(3))=h^0(S,C+4K_S)=1$ by \eqref{eqn:obstruction space}.
Since $-3K_S-C$ is not effective, 
we have $h^0(\mathbb P^3,\mathcal I_C(3))=1$.
Thus \eqref{item:dimension2} follows from \eqref{item:dimension}.
\qed

\medskip

The following example was studied in detail in \cite{Kleppe83}
(see also \cite{Kleppe87}).

\begin{ex}[Kleppe]
  \label{ex:Kleppe}
  Let $S$ be a smooth cubic surface, 
  $E_1$ and $E_2$ two skew lines on $S$
  and $C$ a smooth connected curve on $S$ such that
  $C \sim -4K_S+2E_1+2E_2$, i.e.,~$C \sim (12;4,4,4,4,2,2)$.
  We see that $C$ is of degree $d=16$ and genus $g=29$.
  Since $g<3d-18$, the $3$-maximal family $W:=W(12;4,4,4,4,2,2)$
  containing $C$ is not a component of $(H(16,29)^{sc})_{\red}$
  (cf.~Remark~\ref{rmk:dimension reason}).
  $H(16,29)^{sc}$ has a singularity of codimension $1$ along $W$.
  In fact, we see that 
  $h^1(\mathbb P^3,\mathcal I_C(3))=2$ and $h^1(C,\mathcal O_C(3))=1$.
  Then by proposition~\ref{prop:codimension},
  we have $\dim_{[C]} H(16,29)^{sc}=64$.
  Here this number $64$ equals to
  the expected dimension $4d$ of $H(16,29)^{sc}$ at $[C]$.
\end{ex}

\bibliography{mybib}
\bibliographystyle{abbrv}

\end{document}